%% file: bandit_many_arms.tex
\documentclass[11pt,preprint]{imsart}

\input{./style.tex}

\alglanguage{pseudocode}

\begin{document}

\begin{frontmatter}

	\title{Adaptive Algorithms for Infinitely Many-Armed Bandits: A Unified Framework}
	\runtitle{Infinitely Many-Armed Bandits}
	\begin{aug}
	\begin{aug}
\author{\fnms{Emmanuel} \snm{Pilliat}\ead[label=e1]{emmanuel.pilliat@ensai.fr}}
\address{Univ Rennes, Ensai, CNRS, CREST—UMR 9194
\printead{e1}}
\runauthor{Pilliat}
\end{aug}
	
	\end{aug}

\begin{abstract}
We consider a bandit problem where the budget is smaller than the number of arms, which may be infinite. In this regime, the usual objective in the literature is to minimize simple regret. To analyze broad classes of distributions with potentially unbounded support, where simple regret may not be well-defined, we take a slightly different approach and seek to maximize the expected simple reward of the recommended arm, providing anytime guarantees.
To that end, we introduce a distribution-free algorithm, $\OSE$, that adapts to the distribution of arm means and matches the known minimax lower bounds, up to polylogarithmic factors, in the classical bounded settings, while applying to several broader distribution classes. We characterize the sample complexity through the rank-corrected inverse squared gap function. In particular, we recover known upper bounds and transition regimes for $\alpha$ less or greater than $1/2$ when the quantile function is $\lambda_\eta = 1-\eta^{\alpha}$. We additionally identify new transition regimes depending on the noise level relative to $\alpha$, which we conjecture to be nearly optimal.
Additionally, we introduce an enhanced practical version, $\PROSE$, that achieves state-of-the-art empirical performance for the main distribution classes considered in the literature.
\end{abstract}

\end{frontmatter}

\maketitle

\section{Introduction}

In this manuscript, we consider a problem where the learner has a very limited budget to pull arms from a reservoir of large size, potentially infinite. In contrast to classical bandit problems where the budget vastly exceeds the number of arms, there is no hope of finding the best arm as it is impossible to explore all arms. Nevertheless, when not all arms are equivalent and a non-negligible proportion of arms provide significantly better rewards on average, it is possible to detect some of them without exhaustive exploration.
This setting is of considerable practical interest, as infinitely many-armed bandit problems arise naturally when available options vastly exceed feasible trials. Applications include content recommendation, where tens of thousands of new podcasts are released monthly \cite{aziz2022identifying}, hyperparameter optimization over large configuration spaces \cite{li2018hyperband}, and crowdsourcing tasks like The New Yorker Cartoon Caption Contest with thousands of submissions \cite{jain2020new}. Additional domains include labor markets and resource exploration \cite{wang2008algorithms}. The common challenge is identifying high-quality arms from a large set using a limited sampling budget.

Formally, we consider an infinite and countable set of arms $\mathcal{A}$ from which a learner iteratively chooses which arm to pull. We assume that each arm $a \in \mathcal{A}$ produces rewards $(X_{a,s})_{s \geq 1}$ distributed according to
\begin{equation}\label{eq:distrib_arms}
X_{a,s} = \lambda_{\gamma(a)} + \varepsilon_{a,s} \enspace,
\end{equation}
where $(\gamma(a))_{a \geq 1}$ are independent and uniformly distributed in $[0,1]$, and $\eta \mapsto \lambda_{\eta}$ is a nonincreasing and right-continuous function from $(0, 1]$ to $\mathbb{R}$, which corresponds to the quantile function of the arm means distribution.
The noise variables $(\varepsilon_{a,s})_{s \geq 1}$ are assumed to be independent, centered, and $\zeta^2$-sub-Gaussian, i.e.,
$\mathbb{E}[e^{\theta\varepsilon_{a,s}}] \leq \exp(\zeta^2\theta^2/2)$ for any $\theta \in \mathbb{R}$. We say that $\gamma(a)$ is the rank of arm $a$, and that $a_1$ is better than $a_2$ if $\gamma(a_1) \leq \gamma(a_2)$.

The model defined by \Cref{eq:distrib_arms} is equivalent to assuming that each arm has a random mean sampled from an unknown distribution $\mathcal{D}$, referred to as the arm reservoir distribution in \cite{carpentier2015simple}.
Indeed, the function $\eta \mapsto \lambda_{\eta}$ precisely characterizes the distribution $\mathcal{D}$ of the arm means in $\mathcal{A}$. Starting from any distribution $\mathcal{D}$ on $\mathbb{R}$, it is always possible to define a nonincreasing and right-continuous function $\lambda$ such that $\lambda_{\gamma(a)}$ has distribution $\mathcal{D}$. If $F$ is the cumulative distribution function (CDF) of $\mathcal{D}$, then the function $\lambda_{\eta} = \inf\{x \in \mathbb{R}: F(x) \geq 1-\eta\}$ is nonincreasing and right-continuous, and $\eta \mapsto \lambda_{1-\eta}$ is the inverse CDF of $F$. Since $\gamma(a)$ and $1-\gamma(a) \sim \mathcal{U}(0,1)$ are uniform random variables in $[0,1]$,
\begin{equation*}
\mathbb{P}(\lambda_{\gamma(a)}\leq x)=\mathbb{P}(\lambda_{1-\gamma(a)}\leq x)=\mathbb{P}(1-\gamma(a) \leq F(x)) = F(x) \enspace ,
\end{equation*}
which proves that $(\lambda_{\gamma(a)})_{a \in \mathcal{A}}$ are distributed according to $\mathcal{D}$. Notice that we do not assume $\mathcal{D}$ is bounded, since $\lambda_{\eta} \to +\infty$ as $\eta \to 0$ is allowed. From this observation, we say that arm $a$ is a top-$\rho$ arm if its rank satisfies $\gamma(a) \leq \rho$.

At each timestep $t$, the learner selects an arm $\hat a_t \in \mathcal{A}$ to pull, with this decision informed by previous observations. Upon pulling the chosen arm, the learner receives the observation
\begin{equation}\label{eq:observation_model}
X_t = \lambda_{\gamma(\hat a_t)} + \varepsilon_t \enspace .
\end{equation}
Equivalently, we can express this using the counting variable $N_{a,t} = \sum_{t' = 1}^t \mathbf{1} \{ \hat a_{t'} = a\}$, which tracks how many times arm $a$ has been pulled up to time $t$. With this notation, we can write the observation as $X_t = X_{a,N_{a,t}}$, connecting to the formulation in \Cref{eq:distrib_arms}.

In addition to pulling an arm, the learner must provide a \textit{recommendation} $\hat r_t \in \mathcal{A}$.
Since $\mathcal{A}$ is infinite, identifying a best arm satisfying $\gamma(a) = 0$ is infeasible.
Hence, one of the main objectives in this manuscript is to provide a strategy that maximizes the expected simple recommendation reward with high probability, defined as $\lambda_{\gamma(\hat{r}_t)}$. This quantity is well-defined even when $\mathcal{D}$ is unbounded. In the literature, $\mathcal{D}$ is usually assumed to be bounded, and the main objectives are either to minimize the simple regret $\lambda_0-\lambda_{\gamma(\hat r_t)}$ or to maximize the probability of recommending an $\epsilon$-good arm $\mathbb{P}(\lambda_{\gamma(\hat r_t)} \geq \lambda_0 - \epsilon)$ \cite{bubeck2009pure, katz2020true, carpentier2015simple,zhao2023revisiting,jamieson2013finding}.
These quantities are well-defined only if the distribution $\mathcal{D}$ is bounded, i.e., $\lambda_0 < \infty$.
Note that we do not consider cumulative regret or cumulative reward here since \cite{de2021bandits} showed that it is impossible to adapt to $\mathcal{D}$, even for very simple classes of distributions.

Our goal is to design a strategy that maximizes the expected simple recommendation reward $\lambda_{\gamma(\hat{r}_t)}$ without prior knowledge of the distribution $\mathcal{D}$ or of $\eta \mapsto \lambda_{\eta}$. To handle this lack of information, we require our strategy to be adaptive to the unknown distribution $\mathcal{D}$. We characterize this adaptivity through the ranks of the arms: recommendations $\hat r_t$ should have ranks $\gamma(\hat r_t)$ that are as small as possible.
This leads us to the following formal definition of achievable rank sequences.
\begin{definition}[$\delta$-achievable rank sequence]\label[definition]{def:achievable}
An algorithm producing recommendations $(\hat r_t)$ \textit{achieves} a sequence of ranks $(\eta_t)$ with confidence $1-\delta$ if
\begin{equation}
\mathbb{P}(\gamma(\hat r_t) \leq \eta_t \text{ for all } t \geq 1) \geq 1-\delta \enspace .
\end{equation}
When such an algorithm exists, we say the rank sequence $(\eta_t)$ is $\delta$-achievable.
\end{definition}

Thus, if recommendations $(\hat r_t)$ achieve rank sequence $(\eta_t)$, then each recommended arm $\hat{r}_t$ has mean reward in the top $\eta_t$ fraction of all arms. Since $t$ is not predetermined, this corresponds to an anytime guarantee in the classical bandit literature \cite{lattimore2016refinement,jun2016anytime,degenne2016anytime}.
A sequence $(\eta_t)_{t \geq 1}$ is called a \textit{best $\delta$-achievable sequence} if it is $\delta$-achievable and satisfies $\eta_t \leq \eta'_t$ for all $t \geq 1$ whenever $(\eta'_t)_{t \geq 1}$ is another $\delta$-achievable sequence.
Given a fixed confidence parameter $\delta \in (0,1)$, our primary goal is to derive an upper bound that holds for any best $\delta$-achievable sequence of ranks. This upper bound will directly allow us to derive a lower bound on the expected simple recommendation reward $\lambda_{\gamma(\hat r_t)}$.

\subsection{Related Work}
Significant work has been done under assumptions on the distribution $\mathcal{D}$ of arm means, with a focus on the case where $\mathcal{D}$ is uniform \cite{berry1985bandit,bonald2013two}, or more generally when $\mathcal{D}$ is bounded and $\mathbb{P}(\lambda_0-\lambda_{\gamma(a)} \geq \epsilon) \asymp \epsilon^\beta$ for some unknown $\beta > 0$ \cite{carpentier2015simple,wang2008algorithms}. In the latter case, \cite{carpentier2015simple} established that the minimax rate for simple regret is of order $\tfrac{1}{\sqrt{t}} \lor \tfrac{1}{t^{1/\beta}}$ (up to polylogarithmic factors) under standard Gaussian noise.

Recently, the focus has shifted toward distribution-free methods that make minimal assumptions about the underlying reward distribution $\mathcal{D}$. Rather than assuming parametric forms or specific tail behaviors, these approaches develop algorithms and theoretical guarantees that hold for broad classes of distributions \cite{katz2020true,zhao2023revisiting,li2018hyperband,mason2020finding}. A key technique is the bracketing trick introduced by \cite{katz2020true}, which consists of progressively increasing the size of the arm set under consideration.

Most related to our work are the contributions of \cite{zhao2023revisiting}. They consider a fixed-budget $T$ setting where the number of arms $K$ may be significantly larger than the budget. To derive strong performance guarantees, they combine several existing techniques. Starting with the Sequential Halving (SH) algorithm \cite{karnin2013almost}, they adapt it into an anytime algorithm called DSH using the doubling trick---for general discussions on this trick, see e.g., \cite{besson2018doubling}. This yields an algorithm with good guarantees when $K \leq T$. To adapt to the data-poor regime where $K \gg T$, they incorporate the bracketing trick from \cite{katz2020true}. This leads to their final method $\BSH$ and their main result (Theorem 6): an anytime upper bound on the probability of recommending a suboptimal arm, which holds uniformly over all error levels $\epsilon > 0$.

\subsection{Contributions}
First, we take a different approach from the literature by making no assumptions on the distribution of arm means---its support can be unbounded---or on the noise level. Our main criterion is therefore not simple regret but rather simple recommendation reward, where good performance means recommending arms $a$ with small rank $\gamma(a)$. We characterize the achievable recommendation reward through the sample complexity of obtaining good rank sequences, formalized via a rank-corrected inverse squared gap function in \Cref{eq:inv_squared_gap_func}.

Second, we connect our main result to the existing literature on many-armed bandits. We recover the distribution-dependent upper bound on simple regret from \cite{zhao2023revisiting} up to polylogarithmic factors, which is essentially tight for Polynomial($\alpha$) discrete distributions (characterized by $\lambda_{\eta}= 1-\sum_{i=1}^K (\tfrac{i}{K})^\alpha \mathbf{1}\{i < K\eta \leq i+1\}$), as shown in their Corollary 6.2. For the continuous counterpart---$\mathrm{Beta}(1, 1/\alpha)$ distributions---we recover the minimax optimal bound (up to polylogarithmic factors) from \cite{carpentier2015simple}.

Beyond existing rates, we identify new transition phases for $\mathrm{Beta}(1, 1/\alpha)$ distributions that reveal the dependence on $\alpha$ and noise level $\zeta$ typically hidden in prior work. Specifically, our upper bounds exhibit a transition phase at timesteps around $t \asymp \left(\frac{\zeta^2}{\alpha^2}\right)^{\frac{1}{1-2\alpha}}$, which we interpret as a saturation effect in \Cref{sec:implications}.

Third, we analyze Pareto-$1/\alpha$ distributions, a class of unbounded distributions for which simple regret is undefined. When the noise level satisfies $\zeta > 1$, our upper bound reveals transition regimes: the noise becomes negligible for identifying good arms once $t \gg \zeta^{1/\alpha}$. Before this threshold, when $\alpha < 1/2$, we establish that the expected simple recommendation reward scales at least as $\left(\frac{t}{\zeta^2}\right)^{\frac{\alpha}{1-2\alpha}}$.

Finally, we introduce $\PROSE$, a practical version of our main algorithm $\OSE$. Empirical comparisons with $\BSH$ from \cite{zhao2023revisiting} on Polynomial($\alpha$) instances demonstrate that $\PROSE$ achieves strong simple regret performance with high computational efficiency. Moreover, unlike $\BSH$, $\PROSE$ matches the long-term cumulative regret performance of the classical UCB algorithm.

\subsection{Technical Overview}

From a theoretical perspective, our main technical achievement is to provide a unified framework for many-armed bandits through the ranks of the arms, the quantile function, and the rank-corrected inverse squared gap function.

Specifically, our main result enables us to recover previous results (up to polylogarithmic factors) in terms of simple regret or probability of recommending an $\epsilon$-good arm, and to analyze new settings where the distribution is potentially unbounded. From our main result, we are able to characterize the order of magnitude of optimal simple regret and rank recommendation for several classically studied distributions. It also allows us to present upper bounds in new settings with unbounded distributions, where we conjecture the bounds to be nearly tight. Moreover, our upper bounds provide explicit dependencies on the noise level $\zeta$ and the distribution parameter.

The proof of our main theorem relies mainly on concentration bounds for, on the one hand, the number of top-$\rho$ arms and, on the other hand, uniform concentration bounds on the noise of explored arms. To focus on the order of magnitude related to the distribution of arm means rather than on polylogarithmic factors---which could be of interest for future work---we do not use refined concentration bounds such as the iterated logarithm law for sums of sub-Gaussian variables, but rather rely on standard Hoeffding and Bernstein inequalities.

From a practical perspective, we provide a computationally efficient implementation of $\PROSE$. This addresses a significant technical challenge: a naïve implementation requires quadratic time complexity in the number of timesteps (since we must compute a maximum at each step). To overcome this efficiency bottleneck, we develop optimized code that computes recommendations by iteratively updating a partition of sorted sets and performing cascading operations between these sets, achieving $O(\log^2 t)$ amortized complexity per iteration. Beyond statistical performance, we also compare the computational efficiency of $\PROSE$ with an optimized Julia implementation of the Bracketing Sequential Halving with doubling trick ($\BSH$) from \cite{zhao2023revisiting}. Implementations of $\PROSE$ and $\BSH$ are available on the author's GitHub webpage.

\section{Main Result}

Recommending a good arm amounts to identifying an arm with a small rank $\gamma(a)$. To quantify the difficulty of this task, we introduce the rank-corrected inverse squared gap function $G$ for any pair of ranks $\rho < \nu$:
\begin{equation}\label{eq:inv_squared_gap_func}
    G(\rho, \nu) = \frac{\zeta^2\nu}{\rho(\lambda_{\rho} - \lambda_{\nu})^2} \lor \frac{1}{\rho} \; .
\end{equation}
If $G(\rho, \nu) \leq t/\psi$ for some $\psi > 0$, we say that rank $\rho$ is $\psi$-significant at time $t$ with respect to rank $\nu$. $G(\rho, \nu)$ can be seen as a measure of the sample complexity for detecting a top-$\rho$ arm among other bad arms $a$ of rank of order $\nu$, e.g., $\gamma(a)\in [\nu, 2\nu]$, and it captures three main effects. First, $\zeta^2/(\lambda_{\rho} - \lambda_{\nu})^2$ represents the squared inverse gap between the top $\rho$-quantile and the top $\nu$-quantile of $\mathcal{D}$. Consider two arms $a_1$ and $a_2$ satisfying $\gamma(a_1) \leq \rho$ and $\gamma(a_2) \geq \nu$. The sample complexity for detecting with high probability that $a_1$ is better than $a_2$ is of order $\zeta^2/(\lambda_{\rho} - \lambda_{\nu})^2$, up to logarithmic factors. Second, the scale factor $\frac{\nu}{\rho}$ accounts for potential needle-in-a-haystack effects. Even when the gap between top $\rho$ and $\nu$ quantiles is significant, identifying a top-$\rho$ arm with respect to arms of rank larger than $\nu$ becomes considerably more challenging when $\rho \ll \nu$. This increased difficulty stems from the relative scarcity of top-$\rho$ arms: for every arm $a$ satisfying $\gamma(a) \leq \rho$, there exist on average $\nu/\rho$ arms satisfying $\gamma(a) \in [\nu, 2\nu]$. Finally, the term of order $1/\rho$ captures a fundamental lower bound: regardless of how large $\lambda_{\rho}$ is relative to $\lambda_{\nu}$, we require at least $\Omega(1/\rho)$ samples to observe at least one top-$\rho$ arm.

Now that $G$ characterizes the sample complexity order for detecting top-$\rho$ arms among $\nu$-bad arms, we can define the sample complexity for detecting a top-$\eta$ arm among all other arms.
\begin{definition}[$\psi$-significant rank at timestep $t$]
We define the sample complexity function for detecting a top-$\eta$ arm as the nonincreasing function $S: (0,1]\mapsto \mathbb{R}$:
\begin{equation}\label{eq:def_S}
S(\eta) = \inf_{\rho < \eta} \sup_{\nu \geq \eta}G(\rho, \nu) \enspace .
\end{equation}
Moreover, for $\psi > 0$, we say that rank $\eta$ is $\psi$-significant at time $t$ if $S(\eta) \leq t/\psi$.
\end{definition}
In other words, a rank $\eta$ is $\psi$-significant at time $t$ if there exists $\rho \leq \eta$ such that rank $\rho$ is $\psi$-significant with respect to any other rank $\nu \geq \eta$. This leads us to define $\eta^*_t(\psi)$ as the smallest rank that is $\psi$-significant at time $t$:
\begin{equation}\label{eq:def_eta_star}
\eta^*_t(\psi) = \inf\big\{\eta \in (0, 1) :~ S(\eta) \leq \frac{t}{\psi}\big\} \enspace .
\end{equation} 
Hence, a rank $\eta$ is $\psi$-significant at timestep $t$ if and only if $\eta \geq \eta^*_t(\psi)$.
The following theorem establishes that if we set $\psi := \psi_{t,\delta}$ as a large polylogarithmic factor in $t/\delta$, then any sequence of $\psi_{t,\delta}$-significant ranks is achievable by some algorithm with probability at least $1-\delta$, in the sense of \Cref{def:achievable}.

\begin{theorem}\label{th:main_thm_new}
Fix any $\delta \in (0,1)$ and $t \geq 1$. Assume that $\psi :=\psi_{t,\delta} \geq 2^{30} \log^4(5t/\delta)$, and set the tuning parameter $\beta_t = 6\log(5t/\delta)$. With probability $1-\delta$, for all $t \geq 1$, the recommendation $\hat r_t$ of \Cref{algo:ose} ($\OSE$) defined in \Cref{sec:algos} satisfies $\gamma(\hat r_t) \leq \eta^*_t(\psi)$. This implies in particular that $\lambda_{\gamma(\hat r_t)} \geq \lambda_{\eta^*_t(\psi)}$.
\end{theorem}

The proof of \Cref{th:main_thm_new} is deferred to \Cref{sec:proof}. Hence, the recommendations of $\OSE$ are among the top $\eta^*_t(\psi)$ proportion of all arms with high probability at least $1-\delta$. To better understand the detectability condition characterized by $\eta^*_t(\psi)$, we now present a slightly relaxed version of \Cref{th:main_thm_new}.

\begin{corollary}\label[corollary]{cor:main}
For $\psi:=\psi_{t,\delta} \geq 2^{30} \log^4(5t/\delta)$, we define $\tilde S$ and $\tilde \eta$ as follows:
\begin{equation}\label{eq:relaxed_version}
\tilde S(\eta) = \sup_{\nu \geq 2\eta} G(\eta, \nu) \spaceand \tilde \eta_t(\psi)= 2\inf\{\eta \in (0,\tfrac{1}{2}):~\tilde S(\eta) \leq \tfrac{t}{\psi}\} \enspace .
\end{equation}
Then, $\OSE$ also achieves the sequence $(\tilde \eta_t(\psi))_{t \geq 1}$ with probability at least $1-\delta$.
\end{corollary}

\Cref{cor:main} can be directly deduced from \Cref{th:main_thm_new}. Indeed, $S(2\eta) \leq \tilde S(\eta)$ for any $\eta \in (0, 1/2)$, so that $\eta^*_t \leq \tilde\eta_t$. Hence, with probability at least $1-\delta$, $\gamma(\hat r_t) \leq \eta^*_t \leq \tilde \eta_t$. We refer the reader to \Cref{sec:implications} for direct consequences of \Cref{th:main_thm_new}. These include results for specific distribution classes $\mathcal{D}$: we recover the minimax simple regret upper bound of \cite{carpentier2015simple} (up to polylogarithmic factors) as a special case when $1-\lambda_{\eta} \asymp \eta^\alpha$, and match the more general result of \cite{zhao2023revisiting} on $\epsilon$-good arm identification up to a $\log(1/\delta)$ factor.

In \Cref{sec:implications}, we show that \Cref{th:main_thm_new} recovers the result of \cite{zhao2023revisiting} up to a polylogarithmic factor in $t/\delta$.
Our polylogarithmic term is of order $\log^4(t/\delta)$ and is suboptimal: a more careful analysis could reduce these exponents. However, we accept this loss of polylogarithmic precision in exchange for clearer exposition and simpler analysis. This simplification allows us to highlight the main effect captured by the function $G$ while maintaining practical relevance. Indeed, setting tuning parameters to constants instead of polylogarithmic expressions still achieves state-of-the-art performance in practice with our enhanced algorithm $\PROSE$, described after $\OSE$ in \Cref{sec:algos}.

\section{Algorithms}\label{sec:algos}

Algorithm $\OSE$, which achieves \Cref{th:main_thm_new}, proceeds as follows. At each time $t$, we define an exploration scope---a ``bracket'' in the terminology of \cite{katz2020true, zhao2023revisiting}---consisting of arms eligible for pulling. We pull the arm maximizing the UCB within this scope and recommend the arm with the best LCB overall. 
$\PROSE$ enhances $\OSE$ by ranking arms according to their LCBs and restricting exploration scopes to high-LCB arms. $\PROSE$ yields significantly better empirical performance for a large class of distributions (\Cref{sec:numerical_study}), making it our recommended algorithm in practice.

Let $\mathcal{O}_t \subset \mathcal{A}=\{1, 2,\dots \}$ denote the set of arms observed before time $t$. At each timestep $t$, we either pull an arm from $\mathcal{O}_t$ or observe a new arm $a$. Without loss of generality, we assume $a \leq t$ for any new arm. We denote by $\overline X_{a,t} = \frac{1}{N_{a,t}}\sum_{s=1}^{N_{a,t}}X_{a,s}$ the empirical mean of rewards obtained from arm $a$ up to time $t$. For a given sequence of tuning parameters $\beta_t>0$, we define the upper and lower confidence bounds as
\begin{equation}
\UCB_{a,t} = \overline X_{a,t} + \sqrt{\frac{\zeta^2\beta_t}{N_{a,t}}} \quad \text{ and }\quad
\LCB_{a,t} = \overline X_{a,t} - \sqrt{\frac{\zeta^2\beta_t}{N_{a,t}}} \enspace .
\end{equation}
By convention, we set $\mathrm{UCB}_{a,t} = +\infty$ and $\mathrm{LCB}_{a,t} = -\infty$ for unsampled arms (i.e., when $N_{a,t} = 0$ or, equivalently, when $a > |\mathcal{O}_t|$).
Algorithm \ref{algo:ose} presents our $\OSE$ procedure. At each time step $t > 0$, we sample $U \sim \mathrm{Uniform}[0,1]$ and set $Z = \lfloor t^U \rfloor$ to define the exploration scope size. We pull the arm $\hat{a}_t$ maximizing the UCB within $\{1, \ldots, Z\}$. When $Z > |\mathcal{O}_t|$, this mechanism naturally explores a new arm. We recommend the arm with the highest LCB among all observed arms.

In \Cref{th:main_thm_new}, we set $\beta_t = 6\log(5t/\delta)$ for a given failure probability $\delta >0$. The law of the iterated logarithm for sums of sub-Gaussian random variables (see e.g., Lemma 5 of \cite{verzelen2023optimal}) would allow us to set $\beta_t$ of order $\log(\log(t)/\delta)$, yielding slightly improved polylogarithmic factors as in the BUCB algorithm \cite{katz2020true}. However, we opt for the simpler logarithmic form, which streamlines the proof and aligns with our practical implementation using a constant $\beta_t := \beta$ (see \Cref{sec:numerical_study}). Empirically, for $\PROSE$, the results are not very sensitive to the choice of $\beta$.

\begin{algorithm}[H]
\caption{Optimistic Scope Exploration ($\OSE$) \label{algo:ose}}
\begin{algorithmic}[1]
\Require A set of arms $\mathcal{A}$, a sequence of tuning parameters $\beta_t > 0$
\Ensure A sequence of recommended arms $(\hat r_t)$
\For{$t =1, 2,\dots$}
\State Generate an independent $U\sim \mathcal{U}(0,1)$
\State Set $Z =\lfloor t^U \rfloor$ \Comment{Exploration Scope}
\State Pull arm $\hat a_t = \argmax_{a \leq Z} \UCB_{a,t-1}$ \Comment{Optimistic arm}
\State Recommend $\hat r_t =\argmax_{a \in \mathcal{A}} \LCB_{a,t}$ \Comment{Recommendation}
\EndFor
\end{algorithmic}
\end{algorithm}

$\OSE$ shares strong similarities with Bracketing-UCB ($\BUCB$) from \cite{katz2020true}, which pulls arms maximizing the UCB over independently sampled subsets of exponentially growing sizes. The key difference is that $\OSE$ uses nested subsets $\{1, \dots, Z\}$ with random $Z = \lfloor t^U \rfloor$ for $U$ uniform in $[0,1]$, whereas $\BUCB$ independently samples new subsets at each bracket level with deterministic exponentially growing sizes. We adopt this nested structure because it simplifies the algorithm's presentation and naturally leads to $\PROSE$, which achieves substantially better empirical performance.
Note that $\BSH$ (Bracketing Sequential Halving) from \cite{zhao2023revisiting} uses a similar bracketing trick as $\BUCB$ from \cite{katz2020true} by also iterating over disjoint subsets (brackets) of exponentially growing sizes.
$\OSE$, together with $\BSH$ and $\BUCB$, has two main practical drawbacks. First, they continue pulling arms with small LCBs. Since $Z=1$ occurs with probability at least $1/\log(t)$, approximately $t/\log(t)$ samples are wasted on arm $1$ even when it is clearly suboptimal. Second, the random exploration scope $Z$ has two negative consequences: it increases the variance of $\lambda_{\hat r_t}$ and prevents efficient computation. Computing the argmax of UCBs at each time step results in $O(t^2)$ computational complexity rather than linear.

To address the first drawback, we rerank arms at each step by their LCBs, prioritizing arms identified as promising over randomly chosen ones. For the second drawback, we use deterministic exploration scopes of exponentially growing sizes. In $\PROSE$, we loop over sets of sizes $\lfloor t^{Q(j/\log(t))} \rfloor$, where $Q: [0,1] \to [0,1]$ is a quantile function. Setting $Q(x) = x$ (corresponding to $\mathcal{U}(0,1)$) yields scopes of size $\lfloor e^j \rfloor$, which we recommend for large arm counts. When more aggressive exploration strategies are desired, e.g., if the noise is small so that good arms are easily identifiable, we recommend choosing $Q(x) = x^{1/\gamma}$ with $\gamma > 1$ (the $\mathrm{Beta}(\gamma, 1)$ quantile function). Large $\gamma$ values make $Q(j/\log(t)) \approx 1$, approximating standard UCB by considering nearly all arms at each step.
\begin{figure}[h]
    \centering
    \includegraphics[width=0.5\textwidth]{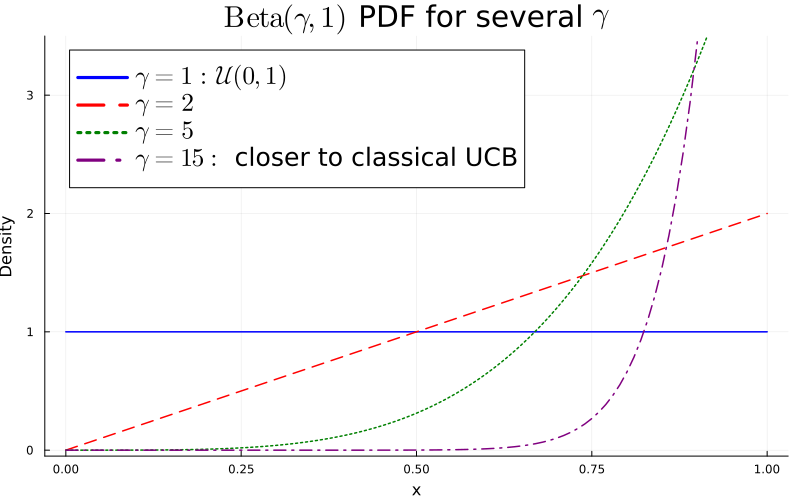}
    \caption{Probability density functions of $\mathrm{Beta}(\gamma, 1)$ distributions for various values of $\gamma$. Higher values of $\gamma$ concentrate the distribution toward 1, corresponding to more aggressive exploration strategies.}
    \label{fig:yourlabel}
\end{figure}
\begin{algorithm}[H]
	\caption{Progressive Ranking for Optimistic Scope Exploration ($\PROSE$) \label{algo:prose}}
	\begin{algorithmic}[1]
		\Require A set of arm $\cA$, a constant tuning parameter $\beta > 0$ to compute the UCBs and LCBs and a function $Q: [0,1] \to [0,1]$.
        \Ensure A sequence of recommended arms $(\hat r_t)$
		\State Initialize $j=1$ and permutation $\pi = \mathrm{id}$ representing the current ranking of the arms
        \For{$t =1, 2,\dots, $}
		\If{$j > \log(t)$}
		\State $j=1$
		\EndIf
		\State Set $Z =t^{Q({j/\log(t)})}$ \Comment{Exploration Scope}
		\State Pull arm $\hat a_t = \argmax_{\pi(a) \leq Z} \UCB_{a,t-1}$ \Comment{Optimistic arm}
		\State Update $\pi$ so that $\LCB_{\pi^{-1}(1),t} \geq \LCB_{\pi^{-1}(2),t} \geq \ldots $ \Comment{Ranking}
		\State Recommend $\hat r_t =\argmax_{a \in \cA} \LCB_{a,t}=\LCB_{\pi^{-1}(1),t}$ \Comment{Recommendation}
		\State $j = j+1$
        \EndFor
	\end{algorithmic}
\end{algorithm}

While \Cref{algo:prose} has a naive computational complexity of $O(t\log t)$ per iteration due to re-sorting all arms by their LCBs in decreasing order and computing the maximum UCB within the exploration scope, we implement in our numerical study a significantly more efficient version with $O(\log^2 t)$ amortized complexity per iteration. The key optimization is to maintain arms in sorted order by LCB in decreasing order using incremental binary search insertions rather than full re-sorting. Additionally, we partition arms into $J \approx \log t$ differential sets: the $j$-th set contains arms whose LCB ranks fall in the interval $(Z_{j-1}, Z_j]$, where $Z_j = \floor{t^{Q(j/\log t)}}$, with each set internally sorted by UCB in decreasing order. Finding the arm with maximum UCB within scope $Z_j$ then requires only $O(J)$ comparisons among the top arms of each set. When pulling an arm changes its confidence bounds and shifts its LCB rank from one interval to another, we perform cascading updates by moving the boundary arms between adjacent differential sets to maintain consistency.

\section{Implications of the Main Theorem}\label{sec:implications}
To fully appreciate the scope of Theorem \ref{th:main_thm_new}, we first demonstrate how it generalizes the main result from \cite{zhao2023revisiting} (up to polylogarithmic factors) before examining specific cases that make distributional tail assumptions, as in \cite{jamieson2013finding,carpentier2015simple}.

\subsection{Uniform $\epsilon$-error probability bound}
To compare with \cite{zhao2023revisiting}, we assume that the distribution of arm means is bounded, which translates to $\lambda_{0}<+\infty$. We define the gap between $\lambda_0$ and $\lambda_{\eta}$ as $\Delta_{\eta} := \lambda_0 - \lambda_{\eta}$. We also define $g(\epsilon) = \sup\{\eta \in [0,1]: \Delta_{\eta} \leq \epsilon \}$, which represents the largest rank $\eta$ whose corresponding gap is at most $\epsilon$. This leads us to define the sample complexity $H(\epsilon)$ as
$$H(\epsilon) = \frac{1}{g(\epsilon/2)} \lor \sup_{\eta \geq g(\epsilon)}\frac{\zeta^2\eta}{g(\epsilon/2)\Delta_{\eta}^2} \; .$$
The supremum of $\tfrac{\zeta^2\eta}{g(\epsilon/2)\Delta_{\eta}^2}$ over all $\eta \geq g(\epsilon)$ is the continuous analogue of the function $H_2$ from \cite{zhao2023revisiting}. The additional term $1/g(\epsilon/2)$ extends our analysis to cases where $\lambda_0 > \zeta$, allowing for arm means that are potentially very large relative to the noise level.
Note that \cite{zhao2023revisiting} restricts their analysis to $\lambda_0 \leq \zeta = 1$, making the term $1/g(\epsilon/2)$ redundant in their setting. Indeed, when $\lambda_0 \leq \zeta = 1$, we have $\tfrac{\zeta^2}{\Delta_1^2} \geq 1$, which implies $\tfrac{1}{g(\epsilon/2)} \leq \tfrac{\zeta^2\eta}{g(\epsilon/2)\Delta_{\eta}^2}$ for $\eta =1$.
While our bound does not explicitly contain the $1/\epsilon^2$ term appearing in the trivial regime of Theorem 6 in \cite{zhao2023revisiting}, this term is implicitly captured by $H(\epsilon)$. Indeed, since $g(\epsilon/2) \leq g(\epsilon)$ and $\Delta_{g(\epsilon)} \leq \epsilon$, we have $H(\epsilon) \geq \tfrac{\zeta^2 g(\epsilon)}{g(\epsilon/2)\Delta_{g(\epsilon)}^2} \geq \zeta^2/\epsilon^2$.

The following corollary establishes that for any $\epsilon>0$, \Cref{algo:ose} returns $\epsilon$-good recommendations $\hat r_t$ with high probability once $t \gg H(\epsilon)$.
\begin{corollary}\label[corollary]{cor:epsilon_gaps}
Fix any $\epsilon > 0$ and $\delta \in (0,1)$. Assume that $\psi \geq 2^{30} \log^4(5t/\delta)$, and set the tuning parameter $\beta = 6\log(5t/\delta)$. With probability at least $1-\delta$, for any $t \geq 4\psi H(\epsilon)$, we have
$$
\lambda_{\gamma(\hat r_t)} \geq \lambda_0-\epsilon \enspace .
$$
\end{corollary}
\Cref{cor:epsilon_gaps} is analogous to Theorem 6 of \cite{zhao2023revisiting}, which analyzes the anytime version of Bracketing Doubling Sequential Halving ($\BSH$). Our result requires a polylogarithmic factor of order $\log^4(t/\delta)$, while \cite{zhao2023revisiting} achieves $\mathrm{polylog}(t) + \log(1/\delta)$. Our proof techniques could yield a $\delta$-dependence of order $\log^2(1/\delta)$, but whether $\OSE$ can match the $\log(1/\delta)$ rate of $\BSH$ remains an open question. Nevertheless, as mentioned previously, we simplify all polylogarithmic factors throughout this manuscript for clearer exposition, and our numerical study demonstrates that these factors are negligible in practice.

The proof of \Cref{cor:epsilon_gaps} relies on a key observation: there exists a gap of at least $\epsilon/2$ between the quantiles corresponding to ranks $g(\epsilon/2)$ and $g(\epsilon)$. When $H(\epsilon) < t/\psi$, we can establish that $S(g(\epsilon)) \leq t/\psi$, which yields $\eta^* \leq g(\epsilon)$.

\begin{proof}[Proof of \Cref{cor:epsilon_gaps}]
Let $\nu \geq g(\epsilon)$ and $\rho < g(\epsilon/2)$. It holds that 
\begin{align*}
G(\rho, \nu) 
= \frac{1}{\rho} \lor \frac{\zeta^2\nu}{\rho(\lambda_{\rho} - \lambda_{\nu})^2} = \frac{1}{\rho} \lor \frac{\zeta^2\nu}{\rho(\Delta_{\nu}-\Delta_{\rho})^2} 
\leq \frac{1}{\rho} \lor \frac{4\zeta^2\nu}{\rho\Delta_{\nu}^2} \; .
\end{align*}
In the last inequality, we used that $\Delta_{\rho} \leq \epsilon/2 \leq \Delta_{\nu}/2$ from the definition of $g$. Hence, taking the limit as $\rho \to g(\epsilon/2)$ and the supremum over $\nu \geq g(\epsilon)$, we obtain that $S(g(\epsilon)) \leq 4H(\epsilon) \leq t/\psi$, which implies that $\eta^*_t(\psi) \leq g(\epsilon)$.
From \Cref{th:main_thm_new}, we obtain that $\gamma(\hat r_t) \leq \eta^*_t(\psi) \leq g(\epsilon)$, and therefore $\lambda_0-\lambda_{\gamma(\hat r_t)} \leq \epsilon$.
\end{proof}

\subsection{Applications to Various Distribution Types}
\Cref{th:main_thm_new} enables us to derive upper bounds on achievable ranks, lower bounds on simple rewards, and, when well-defined, bounds on simple regrets for several distribution classes. In what follows, we analyze three distribution types.

The first is the simplest we can consider: arms have mean $u$ with probability $\eta_0$ or mean $0$ with probability $1-\eta_0$. We refer to this as the Bernoulli-type distribution. This case was precisely analyzed in \cite{de2021bandits} when $\zeta=1$ and $u \in [0,1]$, and we recover their result for identifying a good arm up to polylogarithmic factors in $t/\delta$.

The second class corresponds to $\mathrm{Beta}(1, 1/\alpha)$ distributions: bounded distributions on $[0,1]$ satisfying $\lambda_{\eta}=1-\eta^\alpha$. This class has been extensively studied in the literature, both for bandits with infinitely many arms \cite{carpentier2015simple,bonald2013two} and in its discrete counterpart, the so-called Polynomial($\alpha$) instance problem \cite{zhao2023revisiting,jamieson2013finding}. In extreme value theory, these distributions fall in the Weibull domain of attraction \cite{de2006extreme}. A particularly interesting new result is the phase transition that occurs in our upper bound when $\alpha \geq 1/2$ for Beta distributions, for which we provide intuition below \Cref{eq:ineq_rank_beta}.

The third class consists of Pareto distributions with unbounded support in $[1, +\infty)$, which fall into the Fréchet domain of attraction in extreme value theory. These distributions are characterized by quantile functions $\lambda_{\eta}= \eta^{-\alpha}$ for $\alpha > 0$, and simple regret is not well-defined for them. To the best of our knowledge, no previous results exist for this class, as boundedness of the arm mean distribution is nearly always assumed in these settings. Nevertheless, we identify an interesting phase transition between the regimes $\alpha \leq 1/2$ and $\alpha > 1/2$ when the noise level $\zeta$ is large. The cases $\alpha \leq 1/2$ and $\alpha > 1/2$ respectively correspond to Pareto distributions with lighter and heavier tails.

For both regimes $\alpha \leq 1/2$ and $\alpha > 1/2$ in the Pareto class, we establish an upper bound of order $1/t$ (up to polylogarithmic factors in $t/\delta$) for the rank $\gamma(\hat r_t)$ of the recommended arm when $t \gg \zeta^{1/\alpha}$ (up to a polylog factor). This bound is essentially tight: achieving a better rank than $1/t$ with high probability is impossible, as the probability that no arm $a \leq t$ satisfies $\gamma(a) \leq 1/t$ is at least $1- (1-1/t)^t \geq 1-e^{-1}$.

We call this the quasi-noiseless regime: once $t \gg \zeta^{1/\alpha}$, newly explored arms that are significantly better than the current recommendation can be detected in polylog$(t/\delta)$ steps. While \Cref{th:main_thm_new} yields a trivial upper bound when $\alpha \geq 1/2$ and $t \lesssim \zeta^{1/\alpha}$, it provides a non-trivial bound of order $(\tfrac{\zeta^2}{t})^{1/(1-2\alpha)}$ when $\alpha < 1/2$ and $t \lesssim \zeta^{1/\alpha}$. Although we do not provide optimality results for the Pareto class in this manuscript, this suggests that the transition to the quasi-noiseless regime is much more abrupt when $\alpha \geq 1/2$. The intuition is that the learner must wait time of order $\zeta^{1/\alpha}$ before the distribution $\cD$ outputs arms with mean $\gtrsim \zeta$. We summarize our findings for these three distribution classes in \Cref{tab:sample_table}, with polylogarithmic factors suppressed through the rescaled time $\bt$. 

\begin{figure}[htbp]
    \centering
    \includegraphics[width=\textwidth]{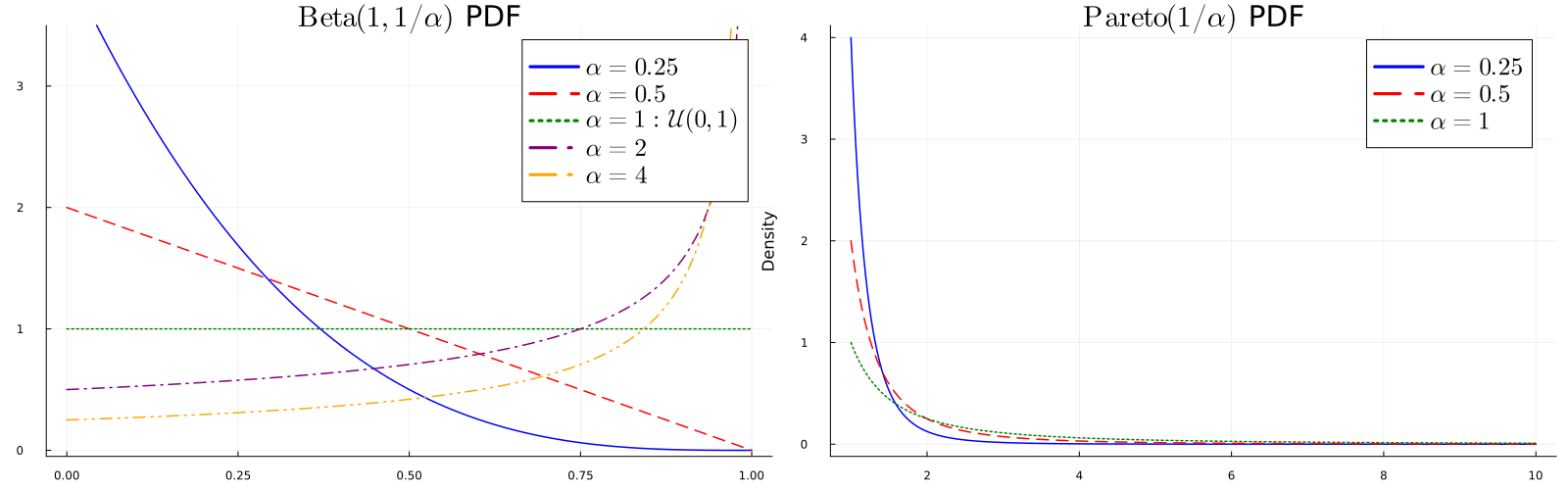}
    \caption{Probability density functions for Beta$(1, 1/\alpha)$ (left) and Pareto$(1/\alpha, 1)$ (right) distributions for various values of $\alpha$. The Beta distribution is supported on $[0,1]$ and becomes more concentrated near 1 as $\alpha$ increases. The Pareto distribution is supported on $[1, \infty)$ and has heavier tail as $\alpha$ increases.}
    \label{fig:beta_pareto_pdfs}
\end{figure}

\begin{table}[h!]
	\centering
	\begin{tabular}{|c c c c|}
	\hline
	 Type of the distribution $\cD$ & Quantile Function $\lambda_{\eta}$ & Upper Bound on $\gamma(\hat r_t)$ & Lower Bound on $\lambda_{\gamma(\hat r_t)}$ \\ \hline

	 \bf{Bernoulli} & $u\1\{\eta \leq \eta_0\}$ & $\eta_0 \1\{\bt \geq \frac{\zeta^2}{\eta_0u^2}\}$ & $u \1\{\bt \geq \frac{\zeta^2}{\eta_0u^2}\}$\\
	
	\bf{Beta} ($\alpha < 1/2$) & $1-\eta^{\alpha}$ & $\frac{1 \lor \zeta^2}{\bt} \lor \left(\frac{\zeta^2}{\alpha^2\bt}\right)^{\tfrac{1}{2\alpha}}$ & $1-\left(\frac{1 \lor \zeta^2}{\bt} \lor \left(\frac{\zeta^2}{\alpha^2\bt}\right)^{\tfrac{1}{2\alpha}}\right)^\alpha$\\
	\bf{Beta} ($\alpha \geq 1/2$) & $1-\eta^{\alpha}$ &$\frac{1}{\bt} \lor \left(\frac{\zeta^2}{\alpha^2\bt}\right)^{\tfrac{1}{2\alpha}} $& $1-\frac{1}{\bt^\alpha} \lor \sqrt{\frac{\zeta^2}{\alpha^2\bt}}$ \\

	\bf{Pareto} ($\alpha < 1/2$) & $\eta^{-\alpha}$ &$\frac{1}{\bt} \lor \left( \frac{\zeta^2}{\bt}\right)^{\frac{1}{1-2\alpha}} $ & $\bt^\alpha \land \left(\frac{\bt}{\zeta^2}\right)^{\frac{\alpha}{1-2\alpha}}$ \\

	\bf{Pareto} ($\alpha \geq 1/2$) & $\eta^{-\alpha}$ & $\frac{1}{\bt}\lor \1\{\bt < \zeta^{1/\alpha}\}$ & $\bt^\alpha \1\{\bt \geq \zeta^{1/\alpha}\}$ \\ \hline
	\end{tabular}
	\caption{Orders of magnitude for the rank $\gamma(\hat r_t)$ and expected simple reward $\lambda_{\gamma(\hat r_t)}$ achieved by $\OSE$ with high probability across different distribution types. Here, $\bt$ denotes $t/ \tilde \psi$, where $\tilde \psi$ is a large polylogarithmic factor in $t/ \delta$ that does not depend on $\alpha$.}
	\label{tab:sample_table}
	\end{table}

\paragraph*{Bernoulli.} 
Assume that $\lambda_{\eta} = u$ if $\eta \leq \eta_0$ and $\lambda_{\eta} = 0$ otherwise. This corresponds to assuming that $\lambda_{\gamma(a)}$ follows a Bernoulli distribution on $\{0,u\}$ with unknown parameter $\eta_0$. When the arm distribution is supported on $[0,1]$ (implying $\zeta \leq 1$), \cite{de2021bandits} established that the best arm can be identified after time $t$ such that $t \geq\log(1/\delta)\log(t)\tfrac{1}{\eta_0 u^2}$. However, they also showed that adapting to the unknown $\eta_0$ is impossible for obtaining non-trivial upper bounds on the cumulative regret, even in this simplest setting.

\Cref{th:main_thm_new} recovers the result of \cite{de2021bandits} up to a polylogarithmic factor. In this setting, $G(\rho, \nu) \leq \tfrac{\zeta^2\nu}{\rho u^2}$ if $\rho \leq \eta_0 < \nu$, and $G(\rho, \nu)=+\infty$ otherwise. Consequently, $S(\eta) = +\infty$ if $\eta \leq \eta_0$, and $S(\eta) = \tfrac{\zeta^2}{\eta_0u^2}$ otherwise. Thus, by \Cref{th:main_thm_new}, as soon as $t \geq \frac{\psi\zeta^2}{\eta_0 u^2}$ -- where $\psi := \psi_{t,\delta}$ denotes the polylogarithmic factor of order $\log^4(t/\delta)$ from \Cref{th:main_thm_new} ---the recommendation $\hat r_t$ satisfies $\gamma(\hat r_t) \leq \eta_0$ with high probability. In other words, one of the good arms among the top $\eta_0$ proportion is identified.

\paragraph*{Beta distribution $(1, \tfrac{1}{\alpha})$.} Let us fix $\alpha > 0$ and assume that $\lambda_{\eta} = 1-\eta^{\alpha}$ for $\eta \in [0, 1]$.
This corresponds to distributions on $[0,1]$ with CDF $F(\lambda) = 1-(1-\lambda)^{1/\alpha}$. Up to constant factors, these distributions correspond to those studied in \cite{carpentier2015simple} with $\alpha = 1/\beta$. The discrete counterpart is called Polynomial($\alpha$) instances in \cite{jamieson2013finding,zhao2023revisiting}.
Then, for $\nu \geq 2\eta$, it holds that

\begin{align*}
	G(\eta, \nu)=\frac{1}{\eta}\lor\frac{\zeta^2\nu}{\eta\left(\nu^{\alpha} - \eta^{\alpha}\right)^2} \leq \frac{1}{\eta}\lor\frac{\zeta^2\nu^{1-2\alpha}}{\eta (1-2^{-\alpha})^2}
\end{align*}

Since the map $\nu \mapsto \nu/(\nu^{\alpha}-\eta^{\alpha})^2$ is decreasing then increasing on $[2\eta,1]$ (its unique interior critical point is a minimum), the supremum $\tilde S(\eta)=\sup_{\nu \geq 2\eta}G(\eta, \nu)$ is attained at an endpoint, $\nu = 2\eta$ or $\nu = 1$. For $\alpha \geq 1/2$ the endpoint $\nu = 2\eta$ dominates; for $\alpha < 1/2$ both endpoints contribute, the one at $\nu = 1$ using the gap $\nu^{\alpha}-\eta^{\alpha}=1-\eta^{\alpha} \geq 1/2$ (valid for $\eta \leq 2^{-1/\alpha}$). Hence, for $\eta \leq 2^{-1/\alpha}$,
\begin{equation}\label{eq:ub_tilde_S}
	\tilde S(\eta)\leq\begin{cases}
		\frac{1\lor 4\zeta^2}{\eta}\lor \frac{8\zeta^2}{\alpha^2 \eta^{2\alpha}}\quad &\text{if } \alpha < 1/2\\
		\frac{1}{\eta}\lor \frac{8\zeta^2}{\alpha^2\eta^{2\alpha}} \quad &\text{otherwise}\\
	\end{cases}
\end{equation}
We used that $1-\eta^{\alpha} \geq 1/2$ for $\eta \leq 2^{-1/\alpha}$ (when $\alpha < 1/2$), and that $2^{\alpha}-1 \geq \alpha/2$ for all $\alpha \in (0,1)$.
Hence, recalling that $\tilde \eta_t(\psi)=2\inf\{\eta \in (0,\tfrac{1}{2}):~\tilde S(\eta) \leq \tfrac{t}{\psi}\}$, we obtain that 
\begin{equation}\label{eq:ineq_rank_beta}
\tilde \eta_t(\psi) \leq
\left\{
\begin{aligned}
&\frac{2\psi(1\lor 4\zeta^2)}{t}\lor \left(\frac{8\psi \zeta^2}{\alpha^2t}\right)^{\frac{1}{2\alpha}} &&\leq \frac{1\lor \zeta^2}{\bt} \lor \left(\frac{\zeta^2}{\alpha^2\bt}\right)^{\frac{1}{2\alpha}} &&\text{if } \alpha < 1/2\\
&\frac{2\psi}{t}\lor \left(\frac{8\psi \zeta^2}{\alpha^2 t}\right)^{\frac{1}{2\alpha}} &&\leq \frac{1}{\bt}\lor \left(\frac{\zeta^2}{\alpha^2\bt}\right)^{\frac{1}{2\alpha}} &&\text{otherwise,}
\end{aligned}
\right.
\end{equation}
where we  set $\bt = t/(24\psi)$.
Applying the relaxed version \Cref{cor:main} of \Cref{th:main_thm_new}, with probability larger than $1-\delta$, we have $\gamma(\hat r_t) \leq \tilde \eta_t(\psi)$.
The approximate inequalities are given up to polylogarithmic factors in $t/\delta$, but the dependence on $\alpha$ captures the correct order of magnitude.

The interpretation is as follows. For $\alpha < 1/2$, the bound is the maximum of a noise-dependent term $\left(\tfrac{\zeta^2}{\alpha^2t}\right)^{1/(2\alpha)}$ and a term $\tfrac{1\lor\zeta^2}{t}$ that depends on the noise only through a constant factor. Since the exponent satisfies $1/(2\alpha) > 1$, the noise-dependent term decreases faster in $t$: it dominates for small $t$ and is overtaken by $\tfrac{1\lor\zeta^2}{t}$ once $t \gtrsim t^* = \left(\tfrac{\zeta^2}{\alpha^2}\right)^{1/(1-2\alpha)}$. This is the mirror image of the case $\alpha \geq 1/2$ analyzed below, where $1/(2\alpha) \leq 1$ and the noise-dependent term instead takes over for $t \gtrsim t^*$, producing a saturation effect. The threshold $\alpha = 1/2$ thus acts as a symmetry axis separating an early-time noise penalty (that washes out) from a late-time one (that saturates).

When $\alpha > 1/2$ and $\zeta \gtrsim \alpha$, the dominant term is of order $\left(\tfrac{\zeta^2}{\alpha^2 t}\right)^{\frac{1}{2\alpha}}$, which corresponds to the regime $\beta < 2$ in \cite{carpentier2015simple}. The interesting behavior occurs when $\zeta \lesssim \alpha$. In this case, there is a phase transition around time $t^*=\left(\tfrac{\zeta^2}{\alpha^2}\right)^{\frac{1}{1-2\alpha}}$. For $t \lesssim t^*$, we obtain the quasi-noiseless bound of order $1/t$, while for $t \gtrsim t^*$, the dominant term becomes $t^{-\frac{1}{2\alpha}}$, as in the $\zeta \gtrsim \alpha$ case.

This phenomenon is somewhat surprising: initially (for $t \lesssim t^*$), the algorithm achieves a very favorable rank guarantee of order $1/t$, but after time $t^*$, the upper bound degrades significantly to $t^{-\frac{1}{2\alpha}}$. While we do not provide a matching lower bound, we conjecture that this is the correct order of magnitude for the rank of the recommendation, and we interpret this phenomenon as a \emph{saturation effect}. For Beta$(1, 1/\alpha)$ distributions with small noise levels $\zeta$, the initial observations before $t^*$ tend to be well separated near the boundary at $1$, which facilitates easy detection of new good arms. However, as time progresses and $t$ exceeds $t^*$, the estimated means of good arms accumulate near $1$, making it increasingly difficult to distinguish among them and identify arms with favorable rank.

To compare with existing results, which do not account for multiplicative factors that depend on $\alpha$ or $\zeta$, let us analyze the expected simple recommendation reward and regret. We derive from the upper bound on $\tilde \eta_t$ that
\begin{equation}
	\lambda_{\gamma(\hat r_t)} \geq 1-\tilde \eta_t^\alpha \geq
	\begin{cases}
		1-\left(\frac{2\psi(1\lor 4\zeta^2)}{t}\right)^\alpha \lor \sqrt{\frac{8\psi \zeta^2}{\alpha^2t}}\quad &\text{if } \alpha < 1/2\\
	1-\left(\frac{2\psi}{t}\right)^\alpha \lor \sqrt{\frac{8\psi \zeta^2}{\alpha^2t}} \quad &\text{otherwise.}\\
	\end{cases}
\end{equation}
Hence, letting $\bt = t/(24\psi)$, we obtain the following simple regret upper bound:
\begin{equation}
	1-\lambda_{\gamma(\hat r_t)} \lesssim
	\begin{cases}
		\left(\frac{1\lor \zeta^2}{\bt}\right)^\alpha \lor \sqrt{\frac{\zeta^2}{\alpha^2\bt}} \quad &\text{if } \alpha < 1/2\\
	\frac{1}{\bt^\alpha} \lor \sqrt{\frac{\zeta^2}{\alpha^2\bt}} \quad &\text{otherwise.}\\
	\end{cases}
\end{equation} 
Notably, when $\zeta=1$, we recover the rate established in \cite{carpentier2015simple} for infinitely-armed bandits, which was shown to be minimax optimal up to polylogarithmic factors. Ignoring constant factors depending on $\alpha$ or $\zeta$, these results also correspond to polynomial($\alpha$) instances in the finite-armed case, as studied in \cite{jamieson2013finding} and \cite{zhao2023revisiting}.

%
%

\paragraph*{Pareto distribution with parameter $\tfrac{1}{\alpha}$.}
Fix $\alpha > 0$ and assume that $\lambda_{\eta} = \eta^{-\alpha}$ for any $\eta \in (0,1)$. This corresponds to unbounded distributions with support on $[1, +\infty)$ and CDF $F(\lambda) = 1-\lambda^{-1/\alpha}$.
For $\nu \geq 2\eta$, we have
\begin{align*}
G(\eta,\nu)=\frac{1}{\eta} \lor \frac{\zeta^2\nu}{\eta\left(\eta^{-\alpha} - \nu^{-\alpha}\right)^2} \enspace .
\end{align*}
The map $\nu \mapsto \nu/(\eta^{-\alpha}-\nu^{-\alpha})^2$ is decreasing then increasing on $[2\eta,1]$ (its unique interior critical point $\nu=\eta(1+2\alpha)^{1/\alpha}$ is a minimum), so $\tilde S(\eta)=\sup_{\nu \geq 2\eta}G(\eta,\nu)$ is attained at the endpoint $\nu=1$. Since $\eta^{-\alpha}-1 \geq \eta^{-\alpha}/2$ for $\eta \leq 2^{-1/\alpha}$, this implies that $\tilde S(\eta) \leq \frac{1}{\eta} \lor 4\zeta^2 \eta^{2\alpha - 1}$.
Hence, $\tilde S(\eta) \leq \tfrac{t}{\psi}$ if $\eta \geq \tfrac{\psi}{t}$ and $4\zeta^2\eta^{2\alpha - 1} \leq \tfrac{t}{\psi}$. When $\alpha \geq 1/2$, a sufficient condition is $\eta \geq \tfrac{\psi}{t}$ and $\tfrac{t}{\psi} \geq 4\zeta^2 (\tfrac{\psi}{t})^{2\alpha - 1}$, which is equivalent to $t \geq \psi (4\zeta^2)^{\frac{1}{2\alpha}}$.

We obtain the following upper bound on $\tilde \eta_t$ depending on whether $\alpha < 1/2$ or $\alpha \geq 1/2$:
\begin{equation}\label{eq:pareto_rank_bound}
\tfrac{1}{2}\tilde \eta_t(\psi) \leq
\left\{\begin{aligned}
&\frac{\psi}{t} \lor \left(\frac{4\psi\zeta^2}{t}\right)^{\frac{1}{1-2\alpha}}
&&\leq \frac{1}{\bt} \lor \left(\frac{\zeta^2}{\bt}\right)^{\frac{1}{1-2\alpha}} &&\text{if } \alpha < 1/2\\
&\frac{\psi}{t}\1\{t \geq \psi (4\zeta^2)^{\frac{1}{2\alpha}}\}
&&\leq \frac{1}{\bt}\1\{\bt \geq (4\zeta^2)^{\frac{1}{2\alpha}}\} 
&&\text{if } \alpha \geq 1/2\enspace ,
\end{aligned}\right.
\end{equation}
where we recall that $\bt = t/ (16 \psi)$.

When $\zeta > 1$ and $\alpha \geq 1/2$, the threshold of order $\zeta^{1/\alpha}$ represents the order of magnitude of the number of arms we must sample before observing an arm with mean at least $\zeta$. Observing an arm with mean exceeding the noise level $\zeta$ before time $t$ of order $\zeta^{1/\alpha}$ is unlikely, as a union bound yields
$$\mathbb P(\exists a \leq t:~ \lambda_{\gamma(a)} \geq \zeta) \leq t\zeta^{-1/\alpha} \enspace .$$
Moreover, before time of order $\zeta^{1/\alpha}$, we have only $\zeta^{1/\alpha} \leq \zeta^2$ samples, which is insufficient to detect a good arm among arms with means in $[0, \zeta]$ under noise level $\zeta$.

However, when $\alpha < 1/2$, the inequality $\zeta^{1/\alpha} \leq \zeta^2$ no longer holds. In this regime, the bound of order $\big(\tfrac{\zeta^2}{\bt}\big)^{1/(1-2\alpha)}$ reflects the fact that even before reaching the time when $\cD$ outputs arms with mean exceeding $\zeta$ (entering a quasi-noiseless regime), it is possible to detect arms with good rank of order $\big(\tfrac{\zeta^2}{\bt}\big)^{1/(1-2\alpha)}$.

Concerning the expected simple reward of the recommended arm, it holds with probability at least $1-\delta$ that
\begin{equation}\label{eq:pareto_reward_bound}
\lambda_{\gamma(\hat r_t)} \geq \begin{cases}
\bt^\alpha \land \left(\frac{\bt}{\zeta^2}\right)^{\frac{\alpha}{1-2\alpha}} \quad &\text{if } \alpha < 1/2\\
\bt^\alpha\1\{\bt \geq (4\zeta^2)^{\frac{1}{2\alpha}}\} \quad &\text{if } \alpha \geq 1/2\enspace ,
\end{cases}
\end{equation}
where $\bt = t/(32\psi)$ (half the value used in the bound on $\tilde \eta_t(\psi)$), and $\psi$ is as defined in \Cref{th:main_thm_new}.

\section{Numerical Study} \label{sec:numerical_study}

In this section, we analyze and compare the performance of four algorithms: $\OSE$ analyzed in \Cref{th:main_thm_new}, its enhanced version $\PROSE$, $\BSH$ from \cite{zhao2023revisiting}, and a naïve implementation of \textbf{UCB} algorithm. We define \textbf{UCB} identically to $\OSE$, except that the exploration scope $Z$ is set to $+\infty$. At each step, it pulls arm with maximum UCB and recommends arm with maximum LCB. For practical simulations, we use constant parameters $\beta$ independent of $t$ and error probability. We recall that our polylogarithmic factors are not optimal and this simplifies tuning to a single parameter. 

We compare these procedures on synthetic data where arm means follow a $\mathrm{Beta}(1, 1/\alpha)$ distribution, corresponding to the quantile function $\lambda_{\eta} = 1-\eta^{\alpha}$. For each arm's reward, we generate normal samples $\mathcal{N}(\lambda_{\gamma(a)}, 1)$. These choices align with standard distributions for arm means \cite{carpentier2015simple,zhao2023revisiting,jamieson2013finding} in the canonical case where $\zeta=1$. 

\begin{figure}[h]
    \centering
    \includegraphics[width=1\textwidth]{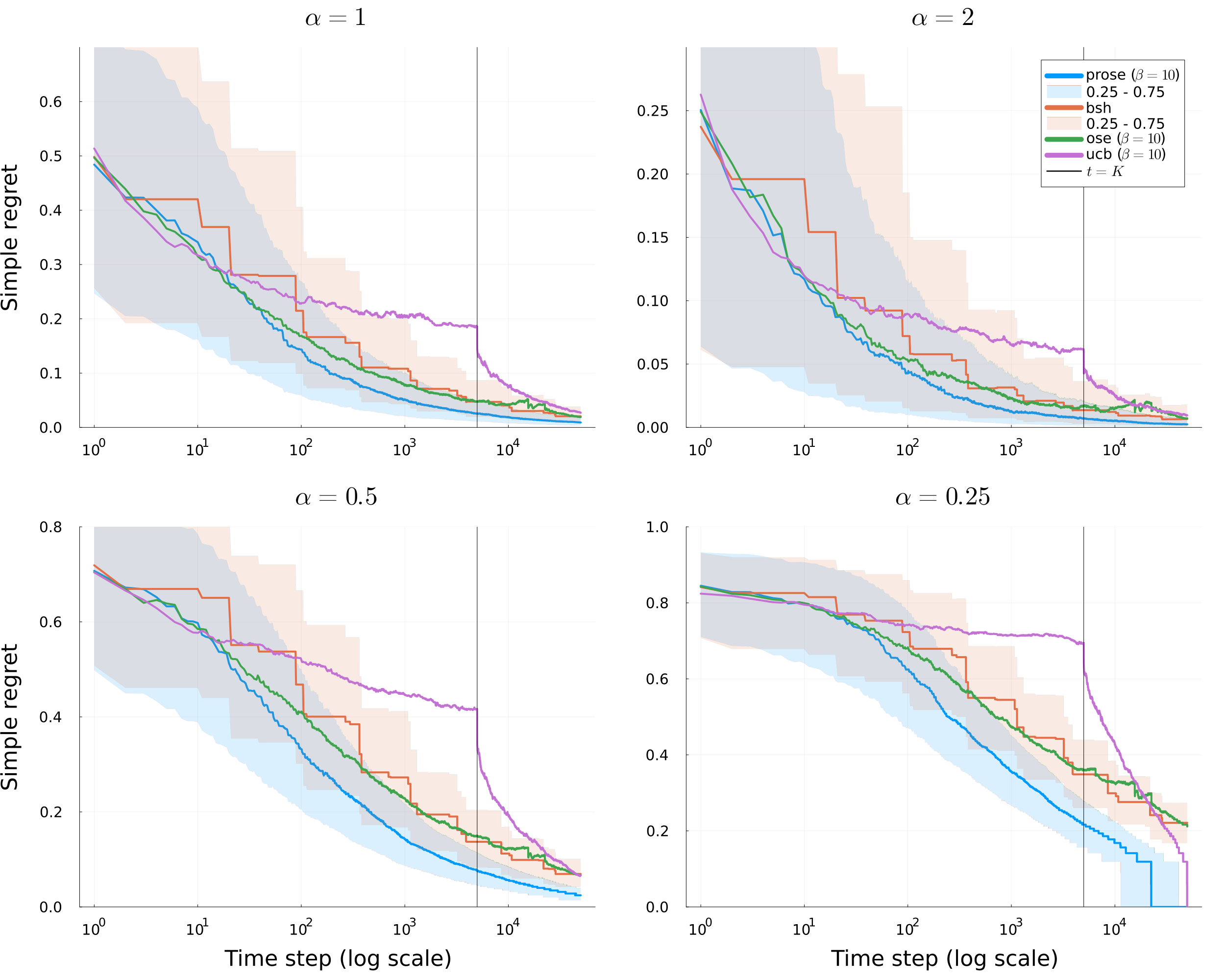}
    \caption{Expected simple regrets for the four algorithms as a function of $t$ on a log scale. Each plot shows results for quantile function $\lambda_{\eta} = 1-\eta^{\alpha}$ with $\alpha \in \{0.25, 0.5, 1, 2\}$, using $K = 5000$ arms and time horizon $T_{\mathrm{max}}=50000$. Trajectories display median simple regret over $N=2000$ trials. Shaded regions (blue for $\PROSE$, orange for $\BSH$) show the interquartile range.}
    \label{fig:simple_regret_main}
\end{figure}

\Cref{fig:simple_regret_main} shows the empirical performance for $\alpha \in \{0.25, 0.5, 1, 2\}$. The results demonstrate that $\PROSE$, with $\beta=10$, consistently outperforms the three other procedures for $t\geq 10$. $\BSH$ exhibits a staircase pattern due to its doubling-trick mechanism. $\OSE$ performs slightly better than $\BSH$ except at large time steps, where $\OSE$ becomes more unstable. As for the naïve UCB algorithm, it only starts to outperform $\OSE$ and $\BSH$ after time $t \gg K$. This behavior is intuitive: before time $t = K$, UCB algorithm performs only pure exploration and starts exploiting better arms only after time $t = K$.

\begin{figure}[h]
    \centering
    \includegraphics[width=1\textwidth]{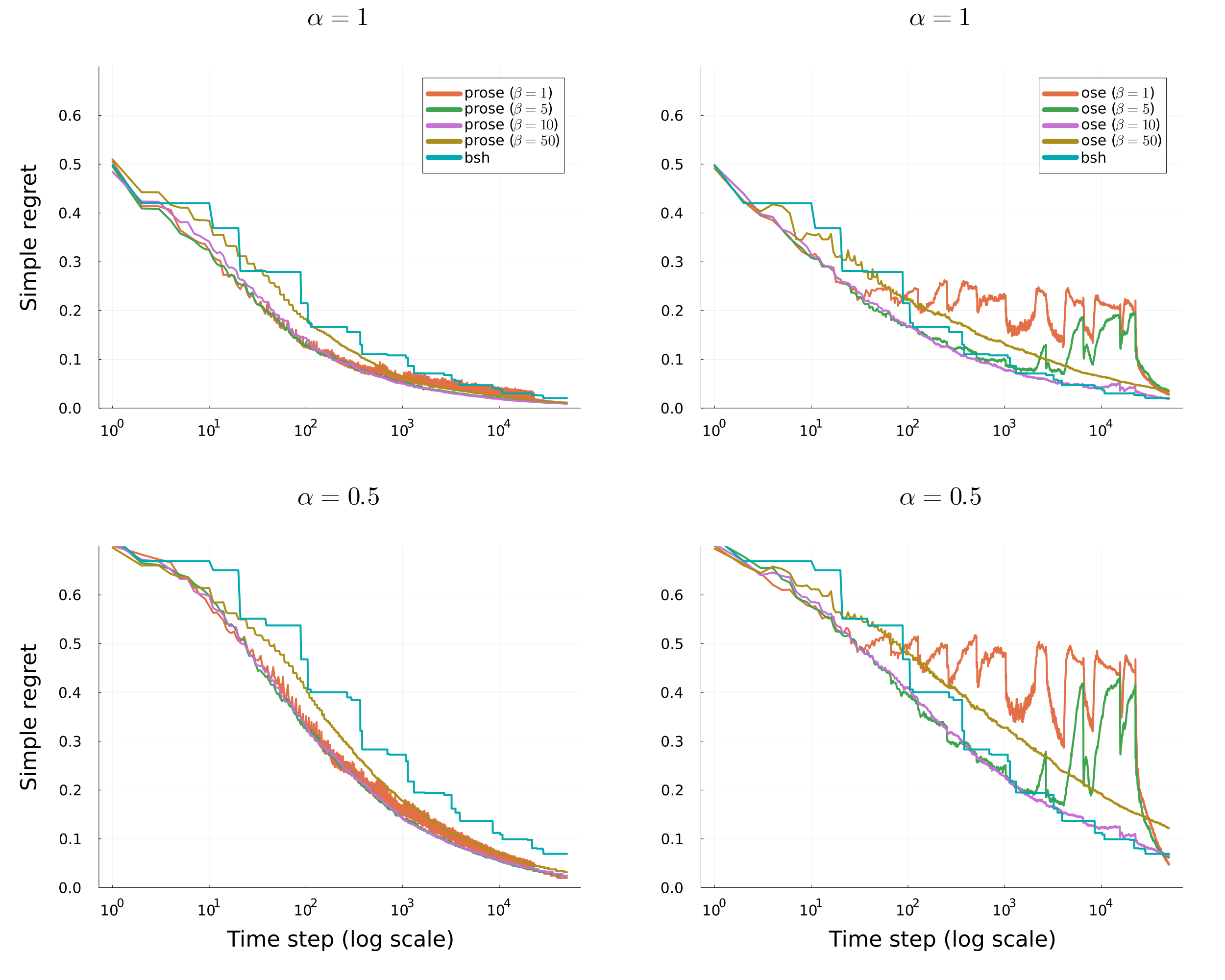}
    \caption{Expected simple regrets for $\BSH$, $\PROSE$ (left), and $\OSE$ (right) with tuning parameter $\beta \in \{1, 5, 10, 50\}$, for $\mathrm{Beta}(1, 1/\alpha)$ instances with $\alpha \in \{0.5, 1\}$. Parameters: $K=5000$ arms, time horizon $T_{\mathrm{max}}= 50000$, and $N=2000$ trials per trajectory.}
    \label{fig:influence_beta}
\end{figure}

The instability of $\OSE$ at large time steps is due to using too small a tuning parameter $\beta = 10$. In theory, this parameter should increase as a polylog of $t$. Surprisingly, this instability is much less pronounced for $\PROSE$, as shown in \Cref{fig:influence_beta} for $\alpha \in \{0.5, 1\}$. This figure shows that while $\OSE$ (right plots) is completely unstable for $\beta = 1$ and too conservative for $\beta=50$, $\PROSE$ (left plots) exhibits only slight instability even for small $\beta$ values and still outperforms $\BSH$ for $\beta=50$.

\begin{figure}[h]
    \centering
    \includegraphics[width=1\textwidth]{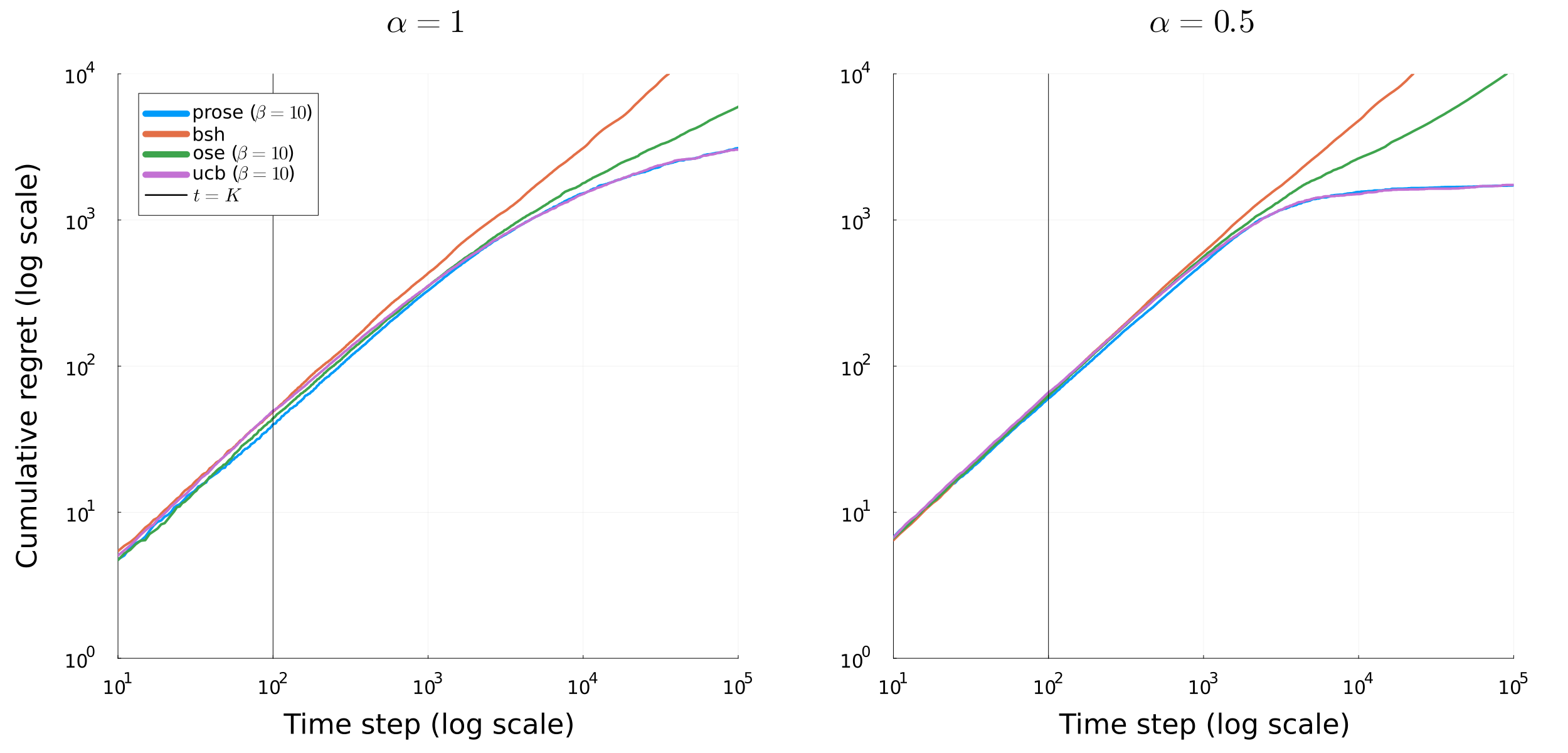}
    \caption{Cumulative regrets for $\BSH$, $\PROSE$, $\OSE$, and \textbf{UCB} on $\mathrm{Beta}(1, 1/\alpha)$ instances with $\alpha \in \{0.5, 1\}$. $\PROSE$ follows the same long-term sublinear trend as \textbf{UCB}, while $\OSE$ and $\BSH$ exhibit quasi-linear growth. Parameters: $K=5000$ arms, $T_{\mathrm{max}}= 50000$, $N=2000$ trials.}
    \label{fig:cumulative_regret}
\end{figure}

While cumulative regret is not a quantity of interest in this paper, it is still instructive to examine how the procedures perform under this metric. Both $\OSE$ and $\BSH$ cannot achieve better than quasi-linear cumulative regret. Indeed, they both continue to use a proportion of order $t/\log(t)$ trials to pull arms that were chosen randomly at the beginning of the learning process.
On the other hand, it is well known---see e.g., Theorem 7.1 of \cite{lattimore2020bandit}---that in the worst case, the \textbf{UCB} algorithm achieves a cumulative regret bound of order $\sqrt{Kt}$ up to polylogarithmic factors, which is sublinear in $t$ when $t \geq K$. We observe this sublinear trend for \textbf{UCB} in \Cref{fig:cumulative_regret}. More interestingly, and perhaps surprisingly, $\PROSE$ follows the same long-term trend as \textbf{UCB}. This demonstrates that, at least for our simulated data, ranking arms at each iteration according to their LCB---the key modification that transforms $\OSE$ into $\PROSE$---achieves best-of-both-worlds performance: strong simple regret when $t \lesssim K$ and recovery of the cumulative regret trend of the \textbf{UCB} algorithm  when $t \gg K$.

All experiments were run on a machine with an Intel Core Ultra 7 165H CPU (16 cores, 32 threads, up to 5.0 GHz) and 31 GB of RAM, running Ubuntu.
For a fair comparison of execution speed, simple regret, and cumulative regret between $\PROSE$ and $\BSH$, we implemented the $\BSH$ algorithm in Julia 1.11 following the additional recommendations from the authors in Appendix D of \cite{zhao2023revisiting}. For $\PROSE$, we use the optimized implementation described at the end of \Cref{sec:numerical_study}, which maintains arms sorted by LCB in decreasing order and partitions them into differential sets sorted by UCB in decreasing order, with all structures kept up to date incrementally. Code is available on the author's GitHub webpage.

Over 2000 trials ($T_{\mathrm{max}}=50000$, $K=5000$ arms), $\PROSE$ averaged 76 ms per trial (range: 30--184 ms) compared to our implemented version of $\BSH$: 89 ms (range: 33--439 ms), representing a 15\% speedup with reduced variance. For comparison, the original $\BSH$ paper \cite{zhao2023revisiting} reported execution times of approximately 2 hours per plot for 500 trials with similar problem sizes ($K \approx 6000$ arms, $T=10000$ timesteps) on an AMD Ryzen 5 PRO 4650GE with 16GB RAM, corresponding to roughly 14 seconds per trial. Our implementation achieves approximately a 100× speedup, likely attributed to a combination of algorithmic optimizations for $\BSH$, modern hardware, and the use of Julia rather than potentially unoptimized implementations in interpreted languages like Python or R.

\clearpage

\bibliographystyle{abbrv}
\bibliography{biblio}

\appendix
\section{Proof}\label{sec:proof}
\subsection{High Probability Events}

In the next three lemmas, we fix $\delta \in (0,1)$, which represents a small failure probability.  Proofs are deferred to \Cref{sec:technical_lemmas}. The first lemma establishes bounds on the noise sums that hold simultaneously for all $t\geq 1$ and all previously observed arms.

\begin{lemma}\label[lemma]{lem:concentration_noise}
	With probability at least $1-\delta/3$, for all $t \geq 1$ and all arm $a \in \mathcal O_t$ observed before time $t$,
	$$ \left|\sum_{s= 1}^{t}\varepsilon_{a,s} \right| \leq \sqrt{6\zeta^2t\log(5t/\delta)}=\sqrt{\zeta^2\beta  t}\; .$$
\end{lemma}

Recall that at each step $t$ of \Cref{algo:ose}, we generate a random variable $Z_t \in [1,t]$ representing the exploration scope, with probability density function $u\mapsto\tfrac{1}{\log(t)}\tfrac{1}{u}$.
For $k\geq 1$, we define $\mathcal T(k, t) = \{t' \leq t:~ Z_{t'} \in [2^{k}, 2^{k+1})\}$, which denotes the set of time steps before $t$ where the exploration scope falls in the interval $[2^k, 2^{k+1}]$.

\begin{lemma}\label[lemma]{lem:control_Z}
	With probability at least $1-\tfrac{\delta}{3}$, for all $t\geq 2^{12}\log^2(t)\log(5t/\delta)$ and all $k \in \{0, \dots, \floor{\log_2(t)}-1\}$,
	$$ |\mathcal T(k, t) | \geq \frac{t}{64\log(t)} \enspace .$$  
\end{lemma}

The following lemma formalizes the intuitive fact that among arms with indices in $[1,z]$ the number belonging to the top $\eta$ fraction is of order $\eta z$ up to log terms and deviations of order $\sqrt{\eta z}$.

\begin{lemma}\label[lemma]{lem:control_quantiles}
	With probability at least $1-\tfrac{\delta}{3}$, for all $t \geq 1$, all $z \in \{1, \dots, t\}$ and all $\eta \in \{2^{-k}, k \in \{0, \dots, \ceil{\log_2(t)}\}\}$,
	$$\left|\sum_{a = 1}^z \1\{\gamma(a)\leq \eta\} - \eta z\right| \leq \sqrt{8\eta z \log\bigg(\frac{5t}{\delta}\bigg)} + 4\log\bigg(\frac{5t}{\delta}\bigg) \enspace .$$  
\end{lemma}
In particular, for any $\eta \in [\tfrac{1}{t}, 1]$, letting $\eta'=2^{-\floor{\log_2(1/\eta)}}$, we have that 
\begin{equation}\label{eq:quantileUB}
	\sum_{a = 1}^z \1\{\gamma(a)\leq \eta\} \leq \sum_{a = 1}^z \1\{\gamma(a)\leq \eta'\} \leq 8\left(\log\bigg(\frac{5t}{\delta}\bigg)\lor 2\eta z\right) \enspace .
\end{equation}
For the second inequality, we used that $1+4+\sqrt{8} \leq 8$ and $\eta' \leq 2\eta$.
Moreover, if $\eta z \geq 64\log(5t/\delta)$, we have with $\eta'' = 2^{-\ceil{\log_2(1/\eta)}}$ that for all $z=1,\dots, t$,
\begin{equation}\label{eq:quantileLB}\sum_{a = 1}^z \1\{\gamma(a)\leq \eta\} \geq \sum_{a = 1}^z \1\{\gamma(a)\leq \eta''\} \geq 1 \enspace .
\end{equation}

\subsection{Analysis of the Procedure, Proof of \Cref{th:main_thm_new}}

We now proceed to the analysis of \Cref{algo:ose}. In what follows, we assume that the three events from \Cref{lem:concentration_noise,lem:control_quantiles,lem:control_Z} simultaneously hold. This occurs with probability at least $1-3\delta/3= 1-\delta$. 
Fix any step $t$ of the procedure and $\eta$ such that $S(\eta) \leq t/\psi$, where $S(\eta)$ is defined in \Cref{eq:def_S}. We aim to show that the recommendation $\hat r_t$ belongs to the top $\eta$ fraction, i.e. $\gamma(\hat r_t) \leq \eta$.

\paragraph{Step 1: Existence of a top-$\rho$ arm with small exploration scope.}

The condition $S(\eta) \leq t/\psi$ implies the existence of $\rho < \eta$ such that $G(\rho, \nu) \leq t/\psi$ for all $\nu \geq \eta$. In other words, there exists $\rho < \eta$ such that for any $\nu \geq \eta$,
\begin{equation}\label{eq:condition_gap}
		\lambda_{\rho} \geq \lambda_{\nu} + \zeta\sqrt{\frac{\nu \psi}{\rho t}} \quad \text{ and } \rho \geq \frac{\psi}{t}\enspace .
\end{equation}
Using \Cref{eq:quantileLB} -- which is derived from \Cref{lem:concentration_noise} -- with $z=\ceil{64\log(5t/\delta)/\rho}\leq t$, there exists an arm $a_{\rho}$ such that
\begin{equation}\label{eq:def_a_rho}
		a_{\rho} \leq 64\log\left(\frac{5t}{\delta}\right)\frac{1}{\rho} \quad \text{ and } \quad \gamma(a_{\rho}) \leq \rho \enspace .
\end{equation}
The insight of \Cref{eq:def_a_rho} is that for an exploration scope of order significantly larger than $1/\rho$, there is at least one arm $a_{\rho}$ in the top $\rho$ fraction.
From the assumption \Cref{eq:condition_gap}, since $\psi \geq 256\log(5t/\delta)$, we have $a_{\rho} \leq t/4$.
Let $k_{\rho} = \floor{\log_2(a_{\rho})}+1$, so that
$$
2^{k_{\rho}-1} \leq a_{\rho} < 2^{k_{\rho}} \enspace .
$$
We recall that $\mathcal T(k_{\rho}, t) = \{t' \leq t:~ Z_{t'} \in [2^{k_{\rho}}, 2^{k_{\rho}+1})\}$. 
We define $N^{(k_{\rho})}_{a,t} = \sum_{t'\in \cT(k_{\rho},t)} \1\{\hat a_{t'} = a\}$, which counts the number of times $t'\leq t$ arm $a$ was pulled when $Z_{t'}\in [2^{k_{\rho}}, 2^{k_{\rho}+1})$. Crucially, whenever $Z_{t'} \in [2^{k_{\rho}}, 2^{k_{\rho}+1})$, the arm $a_{\rho}$ lies within the exploration scope at step $t'$ since $a_{\rho} < 2^{k_{\rho}} \leq Z_{t'}$, and it can thus be pulled if it has the largest UCB at time step $t'$.

	\paragraph{Step 2: Lower and upper bounds on the $\UCB$'s for all arms.}
	Fix any arm $a \geq 1$ and any time step $t' \leq t$. If arm $a$ has not been pulled before time $t'$, then its $\UCB$ at time $t'$ is infinite. Otherwise, $N_{a,t'} \geq 1$ and
	From \Cref{lem:concentration_noise}, we have $|\sum_{s=1}^{N_{a,t'}} \varepsilon_{a,s}| \leq \sqrt{\zeta^2\beta N_{a,t'}}$. Hence,
	\begin{equation*}
		\UCB_{a, t'} = \lambda_{\gamma(a)} + \frac{1}{N_{a,t'}}\sum_{s=1}^{N_{a,t'}} \varepsilon_{a,s} + \sqrt{\frac{\zeta^2\beta}{N_{a,t'}}} \geq \lambda_{\gamma(a)} \enspace .
	\end{equation*}
	In particular, since $\lambda_{\gamma(a_{\rho})} \geq \lambda_{\rho}$,
	\begin{equation}\label{eq:ucb_lb}
		\UCB_{a_{\rho}, t'} \geq \lambda_{\rho} \enspace .
	\end{equation}
	Using the same argument, we also obtain that
	\begin{equation}\label{eq:ucb_ub}
		\UCB_{a,t'} \leq \lambda_{\gamma(a)} + 2\sqrt{\frac{\zeta^2\beta}{N_{a,t'}}} \; .
	\end{equation}
	\paragraph*{Step 3: Non top-$\eta$ arms have not been pulled too often.}
	Let $a$ be an arm that is not in the top $\eta$ fraction, i.e. $\gamma(a)\geq \eta$. 
	Let $t' \leq t$ be the last step in $\cT(k_{\rho},t)$ at which $a$ has been pulled.
	This means that $N^{(k_{\rho})}_{a,t'-1} + 1= N^{(k_{\rho})}_{a,t'} = N^{(k_{\rho})}_{a,t}$. Since $a_{\rho} \leq 2^{k_{\rho}} \leq Z_{t'}$, the $\UCB$ of $a$ at step $t'$ is necessarily larger than the $\UCB$ of $a_{\rho}$ - Otherwise arm $a_{\rho}$ would have been pulled instead:
	\begin{equation*}
		 \UCB_{a,t'-1} \geq \UCB_{a_{\rho},t'-1} \; .
	\end{equation*} 
	Using \Cref{eq:ucb_lb} for $a_{\rho}$ and \Cref{eq:ucb_ub} for $a$ we obtain that $\lambda_{\gamma(a)} + 2\sqrt{\frac{\zeta^2\beta}{N_{a,t'-1}}} \geq \lambda_{\rho}$, which implies that
	\begin{equation}\label{eq:badarm_quickly_dominated}
		N_{a,t}^{(k_{\rho})} = N_{a,t'-1}^{(k_{\rho})} +1 \leq N_{a,t'-1}+1 \leq \frac{4\zeta^2\beta}{(\lambda_{\rho}-\lambda_{\gamma(a)})^2} +1 \leq  \frac{4\beta}{\psi}\frac{\rho t}{\gamma(a)} + 1 \leq \frac{8\beta}{\psi}\frac{\rho t}{\gamma(a)}\; .
	\end{equation}
	In the third inequality, we used assumption \Cref{eq:condition_gap}. In the last inequality, we used that $\tfrac{4\beta\rho t}{\psi \gamma(a)} \geq \tfrac{\rho t}{\psi} \geq 1$.
	Hence, an arm $a$ such that $\gamma(a)\geq \eta$ has been pulled a number of times of order at most $\rho t/\gamma(a)$, within time steps $t'$ in $\cT(k_{\rho}, t)$.
	
	\paragraph*{Step 4: At least one top-$\eta$ arm $a^*$ has been pulled very often.} 

	Since $a_{\rho} \leq t/4$, we have that $k_{\rho} \leq \log_2(t)-1$. Additionally, from \Cref{eq:condition_gap}, $t \geq \psi \geq 2^{12}\log^2(t)\log(5t/\delta)$ and we are in position to apply \Cref{lem:control_Z}:
	\begin{equation}\label{eq:lb_cT}
		|\cT(k_{\rho},t)| \geq \frac{t}{64\log(t)} \enspace .
	\end{equation}
	At any step $t' \in \cT(k_{\rho}, t)$, the arm $\hat a_{t'}$ that was pulled satisfies $\hat a_{t'} \leq Z_{t'} < 2^{k_{\rho}+1}$. Hence,
	\begin{equation*}
		|\cT(k_{\rho},t)|=\sum_{t'\in \cT(k_{\rho},t)}\sum_{a=1}^{2^{k_{\rho}+1}}\1\{\hat a_{t'}=a\} = \sum_{a =1}^{2^{k_{\rho}+1}} N^{(k_{\rho})}_{a,t} \enspace .
	\end{equation*}
	Let $\nu_{l} = 2^{-l}$ for $l \in \cL=\{0, 1, 2, \dots, \floor{\log_2(1/\eta)}\}$, and $\nu_{|\cL|+1}= \eta$. We have

	\begin{equation}\label{eq:ub_cT}
		|\cT(k_{\rho},t)| = \sum_{a =1}^{2^{k_{\rho}+1}} N^{(k_{\rho})}_{a,t}
		= \sum_{l\in \cL}\sum_{a =1}^{2^{k_{\rho}+1}} N^{(k_{\rho})}_{a,t}\1\{\nu_{l+1}<\gamma(a)\leq \nu_{l}\} 
		+ \sum_{a = 1}^{2^{k_{\rho}+1}}N^{(k_{\rho})}_{a,t}\1\{\gamma(a) \leq \eta\}\enspace .
	\end{equation}
	If $\gamma(a) \in (\nu_{l+1}, \nu_l]$, then \Cref{eq:badarm_quickly_dominated} gives that $N_{a,t}^{(k_{\rho})} \leq \frac{8\beta}{\psi}\frac{\rho t}{\nu_{l+1}}$. Hence,
	\begin{align*}
		\sum_{a =1}^{2^{k_{\rho}+1}} N^{(k_{\rho})}_{a,t}\1\{\nu_{l+1}<\gamma(a)\leq \nu_{l}\}
		&\leq \frac{8\beta}{\psi}\frac{\rho t}{\nu_{l+1}}\sum_{a =1}^{2^{k_{\rho}+1}}\1\{\nu_{l+1}<\gamma(a)\leq \nu_{l}\}\\
		&\leq \frac{64\beta}{\psi}(\log(5t/\delta) \lor 2\nu_l 2^{k_{\rho}+1})\frac{\rho t}{\nu_{l+1}}\\
		&\leq \frac{64\beta}{\psi}(\log(5t/\delta) \lor 8\nu_la_{\rho})\frac{\rho t}{\nu_{l+1}}\\ 
		&\leq \frac{2^{13}\beta}{\psi}\log(5t/\delta)t \leq \frac{|\cT(k_{\rho},t)|}{2|\cL|} \; .
	\end{align*}
	In the second inequality, we used \Cref{eq:quantileUB}  with $\nu_l$ and $z=2^{k_{\rho}+1}$. For the fourth inequality, 
	we used that $\nu_l/\nu_{l+1} \leq 2$ and $\rho a_{\rho} \leq 64\log(5t/\delta)$ and that $\rho/\nu_{l+1} \leq 1$ . 
	For the last inequality, we used \Cref{eq:lb_cT} and the fact that $\psi \geq 2^{20}|\cL|\log(t)\log(5t/\delta) \beta $.
	Summing this bound over the $|\cL| \leq \log_2(1/\eta)+1$ rank-shells and plugging into \Cref{eq:ub_cT}, we obtain
	\begin{equation}
		\sum_{a = 1}^{2^{k_{\rho}+1}}N^{(k_{\rho})}_{a,t}\1\{\gamma(a) \leq \eta\} \geq \frac{|\cT(k_{\rho},t)|}{2} \geq \frac{t}{128\log(t)} \enspace ,
	\end{equation}
	where we used \Cref{lem:control_Z} for the last inequality.
	Using \Cref{eq:quantileUB} again, we obtain that
    \begin{align*}
		\frac{t}{128\log(t)} &\leq \left(\max_{a: \gamma(a)\leq \eta} N^{(k_{\rho})}_{a,t}\right)\sum_{a = 1}^{2^{k_{\rho}+1}}\1\{\gamma(a)
		\leq \eta\}\\
		&\leq \left(\max_{a: \gamma(a)\leq \eta} N^{(k_{\rho})}_{a,t}\right)\cdot 8\left(\log\bigg(\frac{5t}{\delta}\bigg)\lor 8\eta a_{\rho}\right) \\
		&\leq 2^{12}\log\bigg(\frac{5t}{\delta}\bigg) \left(\max_{a: \gamma(a)\leq \eta} N^{(k_{\rho})}_{a,t}\right)\frac{\eta}{\rho} \; .
	\end{align*}
	In the last inequality, we used that $a_{\rho} \leq 64\log(\tfrac{5t}{\delta})\tfrac{1}{\rho}$.
	Hence, there exists $a^*$ such that $\gamma(a^*)\leq \eta$ and
	\begin{equation}\label{eq:good_arm_pulled_often}
		N_{a^*,t} \geq N^{(k_{\rho})}_{a^*,t} \geq \frac{1}{2^{19}\log^2(5t/\delta)}\frac{\rho}{\eta}t \enspace .
	\end{equation}
\paragraph{Step 5: This top-$\eta$ arm $a^*$ has an LCB of order $\lambda_{\rho}$.} 

Let $t' \leq t$ be the last step $a^*$ was pulled, so that $\UCB_{a^*, t'-1} = \UCB_{a^*, t-1}$, $\LCB_{a^*, t'} = \LCB_{a^*, t}$ and $N_{a^*, t'-1} + 1 = N_{a^*, t}$.
Then, 
\begin{equation*}
\lambda_{\gamma(a^*)} + 2\sqrt{\frac{\zeta^2\beta}{N_{a^*,t'-1}}} \geq \UCB_{a^*, t'-1} \geq \UCB_{a_{\rho}, t'-1} \geq \lambda_{\rho} \; .
\end{equation*}
Using \Cref{eq:good_arm_pulled_often}, the inequality $\rho \geq \psi/t$, $\eta \leq 1$ and $\psi \geq 2^{20}\log^2(5t/\delta)$, we have that $N_{a^*,t'-1}=N_{a^*, t}-1 \geq \frac{1}{2^{20}\log^2(5t/\delta)}\frac{\rho}{\eta}t$. Hence,
\begin{equation}
\lambda_{\gamma(a^*)} \geq \lambda_{\rho} - \sqrt{2^{20}\log^2(5t/\delta)\zeta^2\beta\frac{\eta}{\rho t}} \enspace .
\end{equation}
This implies that
\begin{equation}
\LCB_{a^*, t} \geq \lambda_{\gamma(a^*)} - 2\sqrt{\frac{\zeta^2\beta}{N_{a^*, t'-1}}} > \lambda_{\rho} - \sqrt{2^{24}\log^2(5t/\delta)\zeta^2\beta\frac{\eta}{\rho t}} \enspace .
\end{equation}

\paragraph{Step 6: All non-top-$\eta$ arms have smaller LCB than $a^*$.}
Let $a$ be any arm in the $(1-\eta)$ bottom fraction, i.e., $\gamma(a) \geq \eta$.
From the same argument as in Step 2, it holds that $\LCB_{a,t} \leq \lambda_{\eta}$.
Hence,
\begin{equation}
\LCB_{a^*,t}- \LCB_{a,t} > \lambda_{\rho}-\lambda_{\eta} - \sqrt{2^{24}\log^2(5t/\delta)\zeta^2\beta\frac{\eta}{\rho t}} \geq \zeta\sqrt{\frac{\eta \psi}{\rho t}} - \sqrt{2^{24}\log^2(5t/\delta)\zeta^2\beta\frac{\eta}{\rho t}} > 0 \enspace .
\end{equation}
The last inequality comes from the fact that we choose $\psi > 2^{24}\log^2(5t/\delta)\beta$ (the factor $\zeta$ cancels on both sides).
Hence, arm $a$ cannot be chosen as recommendation at step $t$ by \Cref{algo:ose}, since $a^*$ has a higher LCB. Therefore, the recommendation $\hat r_t$ is necessarily an $\eta$-good arm, that is $\gamma(\hat r_t) \leq \eta$. As this is valid for any $\eta$ such that $S(\eta) \leq t/\psi$, we conclude that $\gamma(\hat r_t) \leq \eta^*_t(\psi)$.

\section{Proofs of Technical Lemmas}\label{sec:technical_lemmas}
\begin{proof}[Proof of \Cref{lem:concentration_noise}]
	From Hoeffding inequality applied to the independent $\zeta$-subgaussian random variables $\varepsilon_{a,s}$, for any fixed $t \geq 1$, we have with probability at least $1-\tfrac{2\delta}{\pi^2t^3}$ that
	$$ \left|\sum_{s= 1}^{t}\varepsilon_{a,s} \right| \leq \sqrt{2\zeta^2t\log(\tfrac{\pi^2t^3}{2\delta})} \leq \sqrt{6\zeta^2t\log(5t/\delta)} \; .$$
	Since $|\cO_t| \leq t$ and $\sum \tfrac{1}{t^2}=\pi^2/6$, we get the result with a union bound over all possible $a\in \cO_t$ and $t\geq 1$.
\end{proof}
	\begin{proof}[Proof of \Cref{lem:control_Z}]
Let $s \geq 1+3t/4$. Since $2^{k+1} \leq t$, it holds that $s \geq \tfrac{3}{2}2^{k}$. Hence, $Z_s$ satisfies
\begin{equation*}
\mathbb{P}(Z_s \in [2^{k}, 2^{k+1})) = \int_{2^{k}}^{2^{k+1} \land s} \frac{1}{u\log(s)}du = \frac{\log(2^{k+1}\land s) - \log(2^{k})}{\log(s)}\geq \frac{\log(3/2)}{\log(t)} \geq \frac{1}{4\log(t)} \; .
\end{equation*}
The $\floor{t/4}$ random variables $\mathbf{1}\{Z_s \in [2^{k}, 2^{k+1})\}$ for $s \in \{1+\ceil{3t/4}, \dots, t\}$ are independent and take values in $\{0,1\}$. Hence, using Hoeffding's inequality on the bounded variables $\mathbf{1}\{Z_s \in [2^{k}, 2^{k+1})\}$, we obtain that, with probability at least $1-\frac{2\delta}{\pi^2t^3}$,
\begin{equation*}
|\mathcal{T}(k,t)| \geq \sum_{s \geq 3t/4+1}\mathbf{1}\{Z_s \in [2^{k}, 2^{k+1})\} \geq \frac{t}{32\log(t)} - \sqrt{\frac{t}{8}\log\left(\frac{\pi^2t^3}{2\delta}\right)} \geq \frac{t}{64\log(t)}\enspace .
\end{equation*}
For the last inequality, we used the assumption which implies that $\sqrt{\tfrac{t}{8}\log(\tfrac{\pi^2t^3}{2\delta})} \leq \tfrac{t}{64\log(t)}$.
We obtain the result from a union bound over all $t\geq1$ and all $k = 0, \dots, \floor{\log_2(t)}-1$ (there are at most $\log_2(t) \leq t$ such $k$'s).
\end{proof}
	\begin{proof}[Proof of \Cref{lem:control_quantiles}]
Let us fix $z\in \{1, \dots, t\}$ and $\eta \in \{2^{-k}, k \in \{0, \dots, \ceil{\log_2(t)}\}\}$. The random variables $(\gamma(a))_{a=1, \dots, z}$ are iid uniform in $[0,1]$. Hence $\mathbf{1}\{\gamma(a) \leq \eta\}$ are iid Bernoulli random variables with parameter $\eta$.
From Bernstein's inequality---see \cite{boucheron2003concentration}---it holds that, with probability at least $1-\tfrac{2\delta}{\pi^2t^4}$,
\begin{align*}
\left|\sum_{a = 1}^z \mathbf{1}\{\gamma(a)\leq \eta\} - \eta z\right| \leq \sqrt{2\eta z \log\left(\frac{\pi^2t^4}{2\delta}\right)} + \log\left(\frac{\pi^2t^4}{2\delta}\right) \leq \sqrt{8\eta z \log\left(\frac{5t}{\delta}\right)} + 4\log\left(\frac{5t}{\delta}\right)\enspace .
\end{align*}
A union bound over all $z \in \{1, \dots, t\}$ and all $\eta \in \{2^{-k}, k \in \{0, \dots, \ceil{\log_2(t)}\}\}$ gives the same inequality, simultaneously for all $z$ and all $\eta$, with probability $1-\tfrac{2\delta}{\pi^2t^2}$. A final union bound over all $t \geq 1$ gives the result.
\end{proof}

\end{document}

%% file: style.tex
\usepackage{bm}
\usepackage[utf8]{inputenc}
\usepackage{mathrsfs}
\usepackage[T1]{fontenc}

\usepackage{comment}
\usepackage{listings}
\usepackage{algorithmicx}
\usepackage[ruled]{algorithm}
\usepackage{algpseudocode}
\usepackage{algpascal}
\usepackage{algc}

\usepackage{bigints}
\usepackage{amsmath,amsfonts,amssymb}
\usepackage{MnSymbol}

\usepackage{graphicx}
\usepackage[left=3cm,right=3cm,top=2cm,bottom=2cm]{geometry}

\usepackage{hyperref} 
\hypersetup{
    colorlinks=true,
    linkcolor=blue,
    filecolor=magenta,
    urlcolor=cyan,
}

\usepackage{ulem}

\usepackage{amsthm}

\usepackage{cleveref} 
\crefformat{equation}{(#2#1#3)}
\crefformat{line}{#1}
\crefformat{table}{Table #2#1#3}
\usepackage{multirow}
\crefformat{lemma}{Lemma #2#1#3}
\crefformat{definition}{Definition #2#1#3}
\crefformat{proposition}{Proposition #2#1#3}
\crefformat{theorem}{Theorem #2#1#3}
\crefformat{corollary}{Corollary #2#1#3}
\crefformat{question}{Question #2#1#3}

\theoremstyle{definition}
\newtheorem{definition}{Definition}
\theoremstyle{plain}
\newtheorem{theorem}{Theorem}[section]
\newtheorem{corollary}[theorem]{Corollary}

\newtheorem{lemma}[theorem]{Lemma}

\theoremstyle{remark}

\usepackage[english]{babel}
\usepackage[hang,small]{caption}
\usepackage{boxedminipage}
\usepackage{dsfont}

\usepackage[usenames,dvipsnames]{xcolor}
\usepackage{colortbl}

\usepackage[customcolors]{hf-tikz}
\usepackage[textwidth=1.5cm, textsize=scriptsize]{todonotes}
\usepackage{pdfpages}

\newcolumntype{P}[1]{>{\centering\arraybackslash}p{#1}}
\DeclareMathOperator*{\argmax}{arg\, max}

\def\cA{\mathcal{A}}

\def\cD{\mathcal{D}}

\def\cL{\mathcal{L}}

\def\cO{\mathcal{O}}

\def\cT{\mathcal{T}}

\def\bt{\mathbf{t}}

\def\1{{\mathbf 1}}
\def\0{{\mathbf 0}}

\DeclareMathAlphabet\mathbfcal{OMS}{cmsy}{b}{n}

\newcommand{\floor}[1]{\left\lfloor#1\right\rfloor} 
\newcommand{\ceil}[1]{\left\lceil#1\right\rceil}

\def\proscal<#1,#2>{\langle #1,#2\rangle}
\newcommand{\poubelle}[1]{}

\newcommand{\LCB}{\mathrm{LCB}}
\newcommand{\UCB}{\mathrm{UCB}}
\newcommand{\OSE}{\mathbf{OSE}}
\newcommand{\PROSE}{\mathbf{PROSE}}
\newcommand{\BSH}{\mathbf{BSH}}
\newcommand{\BUCB}{\mathbf{BUCB}}
\newcommand{\spaceand}{\quad \text{and} \quad}

\def\bt{\mathbf{t}}

\def\1{{\mathbf 1}}
\def\0{{\mathbf 0}}

\DeclareMathAlphabet\mathbfcal{OMS}{cmsy}{b}{n}